\documentclass[final,5p,twocolumn]{elsarticle}

\journal{System and Control Letters}

\usepackage{amsmath,amssymb}
\usepackage{amsthm}
\usepackage{color}
\usepackage{xcolor}
\usepackage{enumerate}
\usepackage{graphicx}
\usepackage{siunitx}
\usepackage{tikz,pgfplots}
\usepackage[ruled,vlined]{algorithm2e}
\usetikzlibrary{intersections, backgrounds}
\usepgfplotslibrary{groupplots}
\usepackage[utf8]{inputenc}
\usepackage{mathtools}
\usepackage{url}

\newtheorem{lem}{Lemma}
\newtheorem{thm}{Theorem}
\newtheorem{rem}{Remark}
\newtheorem{definition}{Definition}

\makeatletter
\newcommand{\pushright}[1]{\ifmeasuring@#1\else\omit\hfill$\displaystyle#1$\fi\ignorespaces}
\newcommand{\pushleft}[1]{\ifmeasuring@#1\else\omit$\displaystyle#1$\hfill\fi\ignorespaces}
\makeatother

\usepackage{Notation}

\begin{document}

\begin{frontmatter}	
	\title{Structure-Preserving Model Order Reduction for Index One Port-Hamiltonian Descriptor Systems\tnoteref{t1,t2}}
    \tnotetext[t1]{The research by P. Schwerdtner, M. Voigt, and T. Moser was supported by the German Research Foundation (DFG) within the projects 424221635 and 418612884 and that of V. Mehrmann was supported by the DFG through project B03 of SFB TRR 154.}
    \tnotetext[t2]{\textbf{CRediT author statement:} \textbf{Paul Schwerdtner, Tim Moser:} Conceptualization, Methodology, Software, Data Curation, Writing -- Original Draft, Visualization, \textbf{Volker Mehrmann:} Conceptualization, Writing -- Review \& Editing, Supervision, \textbf{Matthias Voigt:} Conceptualization, Writing -- Review \& Editing, Supervision, Funding Acquisition}
    
	\author[TUB]{Paul Schwerdtner\corref{COR}}
	\ead{schwerdt@math.tu-berlin.de}	
	\author[TUM]{Tim Moser\corref{COR}}
	\ead{tim.moser@tum.de}	
	\author[TUB]{Volker Mehrmann}\ead{mehrmann@math.tu-berlin.de}
	\author[FUS]{Matthias Voigt}\ead{matthias.voigt@fernuni.ch}
	
	\affiliation[TUB]{organization={TU Berlin, Institute of Mathematics},
		addressline={Stra{\ss}e des 17.~Juni 136}, 
		city={10623 Berlin},
		country={Germany}}
	\affiliation[TUM]{organization={Technical University of Munich, Chair of Automatic Control},
	addressline={Boltzmannstra{\ss}e 15}, 
	city={85748 Garching/Munich},
	country={Germany}}			
	\affiliation[FUS]{organization={UniDistance Suisse},
	addressline={Schinerstrasse 18}, 
	city={3900 Brig},
	country={Switzerland}}			
	\cortext[COR]{Corresponding authors. Authors contributed equally.}

	\begin{abstract}
          We develop optimization-based structure-preserving model order reduction (MOR) methods for port-Hamiltonian (pH) descriptor systems of differentiation index one. Descriptor systems in pH form permit energy-based modeling and intuitive coupling of physical systems across different physical domains, scales, and accuracies. This makes pH models well-suited building-blocks for component-wise modeling of large system networks. In this context, it is often necessary to preserve the pH structure during MOR\@. We discuss current projection-based and structure-preserving MOR algorithms for pH systems and present a new optimization-based framework for that task. The benefits of our method include a simplified treatment of algebraic constraints and often a higher accuracy of the resulting reduced-order model, which is demonstrated by several numerical examples.

	\end{abstract}
	
\end{frontmatter}

\section{Introduction}%
\label{sec:introduction}
We present optimization-based structure-preserving model order reduction (MOR) algorithms for models described by port-Hamiltonian differential-algebraic equations (pH-DAEs). 
Differential-algebraic equations (DAEs) naturally emerge in the modeling of complex systems because they allow the inclusion of preservation and network laws such as mass-balances in chemical processes, joints in mechanical systems, or Kirchhoff's laws in electrical circuits in the model as algebraic constraints. The use of automatic modeling systems such as \texttt{modelica}\footnote{See \url{https://modelica.org/}.} or \texttt{simscape}\footnote{See \url{https://de.mathworks.com/products/simscape.html}.} has further promoted the use of DAE-based models. 

In recent years, DAE modeling has increasingly addressed the physical properties of the underlying models by incorporating concepts such as passivity or a Hamiltonian structure, leading to pH-DAEs. The concept of pH-DAEs is particularly useful in the modeling of large networks that are constructed from a large number of network components, such as power networks~\cite{MehMOS18}, gas networks~\cite{DomHLMMT21}, or district heating networks~\cite{HauMMMMRS19}. Typically, in such networks the components have widely varying dimensions and different modeling accuracies.
Some models are highly detailed partial differential equation (PDE) systems, typically spatially discretized via finite element, finite difference, or finite volume methods, and other models are surrogate models generated purely from data, see~\cite{MehU22} for a survey of applications.

The port-Hamiltonian paradigm is particularly suited to handling this modeling challenge because it allows for an intuitive energy-based interconnection of systems from different physical domains and of different scale or modeling accuracy, see~\cite{MehU22,JacZ12,SchJ14}. 
A classical example of pH-DAE modeling arises in electrical circuits modeled using modified nodal analysis, \cite{EstT00,GueF99a,GueF99b}.

When the models resulting from the modeling process of complex systems have a large state-space dimension, then the direct simulation or model-based control of such large-scale systems is often infeasible. Then, typically, model order reduction (MOR) is employed to determine an approximation to the given full-order model (FOM) with a smaller state-space dimension that enables efficient simulation and control. However, the need for optimized operation of large networks of complex systems has revealed the need for a more hierarchical modeling approach, see, e.\,g., ~\cite{DomHLMMT21,DaeMRS21,MehSS18}. This approach is often carried over to the MOR of hierarchical systems. Here separate low-order surrogates are computed for the individual subsystems (potentially at different accuracy levels) instead of one reduced-order model (ROM) for the whole system.

The paradigm shift from applying MOR to one (monolithic) system to using MOR to reduce the components of networked models makes the preservation of certain structural properties of the components essential. This is because one network component may rely on the specific properties (such as passivity, see Section~2) of other components. Furthermore, the properties that result from the network structure of these components must be preserved during MOR such that the coupling of the reduced-order components can be performed in the same way as the coupling of their full-order counterparts. The preservation of the pH-DAE structure ensures the preservation of these network-relevant properties and thus enables a hierarchical low-order modeling approach.

However, structure-preserving MOR for pH-DAEs has still only been partially resolved.
MOR methods for pH models based on ordinary differential equations (pH-ODEs), such as 
~\cite{Gugercin2012,Polyuga2012} have been extended to pH-DAEs in~\cite{HauMM19,BeaGM21}, 
but typically the algebraic equations have to be identified and treated separately to prevent destroying the constraint structure,
see~\cite[Remark 3]{BeaGM21}.
An alternative MOR approach for structure-preservation is passivity-preserving MOR (see Section~\ref{sec:preliminaries}). 
However, these methods, such as positive-real balanced truncation (PRBT), as presented in~\cite{Reis2010}, also do not result in a significant reduction in the equations describing the algebraic constraints without further treatment.
A minimal realization of the subsystem corresponding to the algebraic constraints can be determined by solving discrete-time projected Lyapunov equations~\cite{MehS05}, but without preserving the pH structure.
A recently proposed passivity-preserving MOR method for pH-ODEs based on spectral factorization~\cite{BreU21} may overcome this problem but in its current form an extension to DAEs requires system transformations to identify and separately deal with the constraint equations.

We consider linear constant-coefficient pH-DAEs, defined as follows.
\begin{definition}\label{def:pHsys} \textup{\cite{HauMM19,BeaMXZ18}} A linear constant coefficient DAE system of the form
  \begin{align}\label{eq:FOM}
   \begin{split}
      E \dot x(t) &= (J-R)Qx(t) +(G-P)u(t), \\
      y(t) &= {(G+P)}^\T Qx(t) +(S-N)u(t),
   \end{split}
  \end{align}
  where $E,\,Q,\,J,\,R \in \R^{\nx \times \nx}$, $G,\, P \in \R^{\nx\times \enu}$, $S,\,N \in \R^{\enu\times \enu}$, is called a \emph{port-Hamiltonian differential-algebraic equation (pH-DAE)}, if the following conditions are satisfied:
  \begin{enumerate}[(i)]
    \item The matrices $Q^\T J Q$ and $N$ are skew-symmetric.
    \item The passivity matrix
      \begin{align*}
        \label{eq:passivitymatrix}
        W_{\rm P} := \mat{Q^\T R Q & Q^\T P \\ P^\T Q & S}
      \end{align*}
      and the product $Q^\T E$ are symmetric positive semi-definite (denoted as $\geq 0$ in the following).
  \end{enumerate}
The \emph{Hamiltonian (energy-storage) function} $\mathcal{H}: \R^\nx  \rightarrow \R$ is then given by
\begin{align*}
  \mathcal{H}(x) = \frac{1}{2}x^\T Q^\T E x.
\end{align*}
\end{definition}
Structure-preserving MOR is aimed at computing systems of the form
\begin{align*}
  \begin{split}
     E_r \dot x_r(t) &= (J_r-R_r)Q_r x_r(t) +(G_r-P_r)u(t), \\
    y_r(t) &= {(G_r+P_r)}^\T Q_r x_r(t) +(S_r-N_r)u(t),
  \end{split}
\end{align*}
where the system matrices $E_r,\, Q_r,\, J_r,\, R_r \in \R^{r \times r}$, $G_r,\, P_r \in \R^{r\times \enu}$, $S_r,\, N_r \in \R^{\enu\times \enu}$ satisfy the structural constraints given in Definition~\ref{def:pHsys} with $r \ll \nx$.

For such systems we develop optimization-based structure-preserving MOR algorithms that
\begin{enumerate}[(i)]
  \item work with the untransformed FOM matrices and often need no state transformation at all,
  \item ensure that the ROM is a pH-DAE, however, without the need for a preservation of the algebraic equations, and 
  \item provide ROMs with high accuracy both in terms of the $\hinf$ and the $\htwo$ error.
\end{enumerate}
In this paper we only discuss the case of pH systems with \emph{differentiation index one}; see \cite{KunM06,Meh15} for a detailed discussion of different indices. A simple characterization of systems with differentiation index one is that if the columns of the matrix $V_E \in \R^{n \times k}$  span the kernel of $E$, then $[E, (J-R)Q V_E]$ is of full rank. Such a system is often referred to as \emph{impulse-free} in the literature. Cases involving a higher differentiation index are more complex and will be treated in a future paper. Moreover, we assume that the pencil $sE - (J-R)Q$ is \emph{regular} (i.\,e., $\det(sE-(J-R)Q)$ is not identically zero for all $s \in \C$) and that it is   \emph{asymptotically stable}
(i.\,e., all its finite eigenvalues have a negative real part). The same assumptions are also imposed on the reduced pencil $sE_r - (J_r-R_r)Q_r$.

The paper is organized as follows: in the next section, we cover objectives and state-of-the-art methods for structure-preserving MOR\@. In Section~\ref{sec:our_method}, we explain our optimization-based approach for MOR\@. In particular, we extend  previous work~\cite{Moser2020,Schwerdtner2020} to the DAE case.
Finally, the effectiveness of the proposed methods is demonstrated by an number of numerical experiments.

\section{Preliminaries}%
\label{sec:preliminaries}
We focus on linear time-invariant pH-DAEs of form~\eqref{def:pHsys} with $\pQ=I_n$ to simplify the presentation in this section. Note that it has been shown in~\cite{MehMW21} that it is always possible to achieve this simplification. If $\pQ$ has full column rank, then this is achieved by merely multiplying the system by $\pQ^\T$ from the left and then renaming the system matrices, or alternatively, the part associated with the kernel of $\pQ$ can be removed without changing the Hamiltonian, see~\cite{MehU22} for a detailed discussion.

To approximate the input-to-output behavior of the FOM, we make use of the \emph{transfer function} $\pHtf$ which for~\eqref{eq:FOM} (under our assumption that $\pQ=I_n$) is defined as
\begin{equation*}
	\pHtf(s) = (\pG+\pP)^\T(s\pE-(\pJ-\pR))^{-1}(\pG-\pP)+(\pS-\pN).
\end{equation*}
This function is well-defined because we have assumed that the system is regular and hence it is a matrix with real-rational functions as entries.
The transfer function $\pHtfr$ of the ROM is defined analogously.   Any rational transfer function can be decomposed such that 
\begin{equation*}
    H(s) = H_{\rm sp}(s)+H_{\rm pol}(s),
\end{equation*}
where $H_{\rm sp}$ denotes the \emph{strictly proper} part with $\lim_{s\rightarrow \infty} H_{\rm sp}(s) = 0$ and where $H_{\rm pol}$ is a  matrix polynomial, see, e.\,g., \cite{Antoulas2005}. Based on the condition that the uncontrolled DAE is of index one, $H_{\rm pol}(s) \equiv \mathcal{D}_0 \in \R^{m \times m}$ is constant. Since the structural properties of a pH-DAE  and its Hamiltonian are preserved under a change of basis and a scaling of the equation with an invertible matrix, this decomposition of the transfer function may be obtained by transforming the index-one pH-DAE to the \emph{semi-explicit form}
\begin{align*}
  \begin{split}
    \mat{E_{11} & 0 \\ 0 & 0} \mat{\dot x_1(t) \\  \dot x_2(t)} &= \mat{L_{11} & 0 \\ L_{21} & L_{22}} \mat{x_1(t) \\ x_2(t)} + \mat{G_1 - \pP_1 \\ G_2 - \pP_2}u(t), \\
    y(t) &= \mat{{G_1 + \pP_1} \\ {G_2 + \pP_2}}^\T \mat{x_1(t) \\ x_2(t)} + (\pS - \pN) u(t),
  \end{split}
\end{align*}
where
\begin{align*}
\mat{L_{11} & 0 \\ L_{21} & L_{22}} =  \mat{J_{11} & -J_{21}^\T \\ J_{21} & J_{22}}-\mat{R_{11} & R_{21}^\T \\ R_{21} & R_{22}},
\end{align*}
and where $L_{22}$ and $E_{11}$ have full rank (see \cite{BeaMXZ18} for details). Then, the transfer function $H_{\rm sp}$ is the transfer function of the implicit pH-ODE system
\begin{equation} \label{eq:FOM_decpl}
	\begin{aligned}
	  \begin{split}
	  E_{11} \dot x_1(t) &= (J_{11}-R_{11})x_1(t) + (G_{\rm sp} - P_{\rm sp})u(t), \\
	  y_{\rm sp}(t) &= {(G_{\rm sp} + P_{\rm sp})}^\T x_1(t),
	  \end{split}
	\end{aligned}
\end{equation}
where
\begin{align*}
  G_{\rm sp} &:= G_1 - Y, \quad P_{\rm sp} := P_1 - Y, \\
  Y &:= \frac{1}{2} {(J_{21}-R_{21})}^\T{(J_{22}-R_{22})}^{-\T}(G_2+P_2).
\end{align*}
Furthermore, the constant part of the transfer function is
\begin{align}
  \mathcal{D}_0   &:= S-N-{(G_2+P_2)}^\T{(J_{22}-R_{22})}^{-1}(G_2-P_2). \label{eq:D}
\end{align}
Here, the state $x_2$ is uniquely determined by the algebraic constraint (given by the second equation block) which imposes a consistency condition on the initial value, see \cite{BeaMXZ18}. Note that sparsity patterns that are typically present in the full order matrices cannot in general be preserved with this transformation.

Due to the asymptotical stability of $sE-(J-R)Q$ the transfer function $\pHtf$ is an element of the Hardy space $\mathcal{RH}_\infty^{m \times m}$ of all real-rational $m \times m$ matrix-valued functions which are bounded on the imaginary axis and with all poles having a negative real part. This vector space is equipped with 
the norm
\begin{equation*}
	\quad \left\| \pHtf \right\|_{\mathcal{H}_\infty} := \sup_{\omega \in \R} {\|\pHtf(\mathsf{i}\omega)\|}_2.
\end{equation*}
Moreover, if $\pHtf$ is additionally \emph{strictly proper}, then $H$ is additionally in the Hardy space $\mathcal{RH}_{2}^{m \times m}$ 
which is equipped with the norm
\begin{equation}
  \label{eq:htwonorm}
	\left\| \pHtf \right\|_{\mathcal{H}_2} := {\left(\frac{1}{2\pi}\int_{-\infty}^\infty {\|\pHtf(\mathsf{i}\omega)\|}_{\rm F}^2 \mathrm{d}\omega\right)}^{1/2}.
\end{equation}
We refer to \cite{ZhoDG96} for a detailed discussion of these spaces.

The search for reduced-order models that minimize the error $\pHtf-\pHtfr$ with respect to the $\cH_2$ or $\cH_{\infty}$ norm is dominated by two types of restrictions for pH-DAEs as in \eqref{eq:FOM}. On the one hand, the ROM has to respect the impact of the algebraic constraints in \eqref{eq:FOM} on the input-to-output behavior of the original model. On the other hand, we restrict our search to ROMs which have a pH representation in order to retain the pH structural conditions of the FOM.  
Many existing model reduction methods for pH systems are based on 
the strong connection between the existence of a pH representation and \emph{passivity}, see, e.\,g., \cite{BeaMV19} which is given by the positive real lemma.
\begin{thm}\textup{\cite{BeaMV19}}\label{thm:KYP}
	Suppose that the linear time-invariant model 
		\begin{align}\label{eq:FOM_SS}
		 \begin{split}
			\dot x(t) &= Ax(t) +Bu(t), \\
			y(t) &= Cx(t) +Du(t),
		\end{split}
		\end{align}
	is \emph{minimal}, i.\,e., the pair $(A,B)$ is \emph{controllable} and the pair $(A,C)$ is \emph{observable}. Then there exists a positive definite matrix $X\in\R^{n\times n}$ which satisfies the \emph{Kalman-Yakubovich-Popov (KYP)} linear matrix inequality
	\begin{equation}\label{eq:KYP_ineq}
		\mathcal{W}(X)=\mat{-A^\T X - XA & C^\T-XB \\ C-B^\T X & D+D^\T} \geq 0,
	\end{equation}
	if and only if the system is passive. Moreover, an implicit pH-ODE representation of~\eqref{eq:FOM_SS} may be obtained by setting 
	$\pE=X$, and 
	\begin{equation*}
		\begin{array}{ll}
\pJ=\frac{1}{2}(\pA -A^\T),
& \pG = \frac{1}{2}(\pC^\T+B), \\ 
\pR=-\frac{1}{2}(\pA +A^\T), & \pP = \frac{1}{2}(\pC^\T-B), 
\\ \pN = \frac{1}{2}(D^\T-D), &  \pS=\frac{1}{2}(D^\T+D).
		\end{array}
	\end{equation*}
\end{thm}
Generalizations of the connection between passivity and solvability of a generalized KYP inequality for DAE systems also exist, for recent results, see \cite{ReiRV15,ReiV15}.

The connection between passivity and the pH structure enables two different approaches for reducing pH models: either by directly enforcing a pH structure for the ROM or by applying passivity-preserving MOR methods combined with a subsequent transformation of the ROM to pH form. 
\subsection{PH-preserving MOR techniques}%
Traditional methods which directly retain the pH form in the reduction process are based on Galerkin projections, see, e.g., \cite{Antoulas2020, MehU22} for surveys. The original state $x(t)$ is approximated by ${x(t)\approx V x_r}$,
where the columns of ${V\in\R^{n\times r}}$ form a basis for a suitably chosen subspace of dimension $r$. For instance, in the pH-ODE case with $\pE= I_n$, the ROM coefficient matrices are computed by
\begin{equation*}
	\begin{array}{lll}
		\rJ=U^\T\pJ U, & \rG = U^\T\pG, & \rQ = V^\T\pQ V, \\ \rR=U^\T\pR U, & \rP = U^\T\pP, & \rN=\pN, \quad \rS=\pS,
	\end{array}
\end{equation*}
where ${U = QV(V^\T QV)^{-1}}$, which clearly enforces structure preservation. In \cite{Gugercin2012} an adaptation (IRKA-PH) of the well-known iterative rational Krylov algorithm (IRKA) was proposed that iteratively updates $V$ to fulfill a subset of $\cH_2$ optimality conditions via tangential interpolation of the original transfer function. While IRKA does not ensure stability (or even passivity) of the ROM a priori, it leads to (locally) $\cH_2$ optimal models upon convergence. IRKA-PH, on the other hand, preserves the pH structure and thus produces passive ROMs which, however, generally only fulfill a subset of the $\cH_2$ optimality conditions. Consequently, no $\cH_2$ optimality is achieved in general. The matrix $V$ can also be chosen in order to approximate the Dirac structure of the original model, resulting in the effort- and flow-constraint method \cite{Polyuga2012}. The extension of IRKA-PH and Dirac structure approximation to pH-DAEs is addressed in \cite{HauMM19,BeaGM21}. 

\subsection{Passivity-preserving MOR techniques}%
Positive real balanced truncation (PRBT) for ODEs is a well-studied MOR method, see  \cite{GugA04} and the references therein for a survey. The method is based on computing a minimal solution $X_{\min}$ of \eqref{eq:KYP_ineq} (in the sense of the Loewner ordering in symmetric matrices) and a minimal solution $Y_{\min}$ of the \emph{dual KYP inequality} 
\begin{equation*}
    \mat{-AY - YA^\T & B-YC^\T \\ B^\T-C Y & D+D^\T} \geq 0.
\end{equation*}
The solutions $X_{\min}$ and $Y_{\min}$ are then used to transform the system to a positive-real balanced realization where the transformed minimal solution of the KYP inequalities $\hat{X}_{\min}$ and $\hat{Y}_{\min}$ are equal and diagonal, i.\,e.,  $\hat{X}_{\min} = \hat{Y}_{\min} = \text{diag}(\eta_1,\ldots,\eta_n)$. Then the states corresponding to the small positive real characteristic values $\eta_i$ can be truncated. Finally, PRBT admits an a priori error bound in the gap metric that is derived in \cite{GuiO14}.

Another passivity-preserving MOR technique is  achieved via spectral factorization \cite{BreU21}. It initially requires a factorization 
\begin{equation}\label{eq:specfac}
	\mathcal{W}(X)=\mat{L&M}^\T\mat{L&M}
\end{equation} 
with $L\in\R^{k\times n}$, $M\in\R^{k\times m}$ and where $\text{rank}\mat{L&M}$ is as small as possible. If the pair $(A,B)$ is \emph{stabilizable} and $X = X_{\min}$ is the minimal solution of $\mathcal{W}(X) \ge 0$, then the resulting spectral factor 
system $(A,B,L,M)$
may be reduced via traditional (unstructured) MOR techniques, such as IRKA or balanced truncation (BT), and the passive ROM is obtained from the reduced spectral factor. 

For all passivity-preserving reduction methods, the passive ROM may eventually be transformed back to a pH representation by applying Theorem \ref{thm:KYP}.

\section{A new optimization-based approach}%
\label{sec:our_method}

Instead of obtaining the ROM by projection, we propose using optimization techniques to determine the coefficients of a low-order pH system such that its transfer function matches the transfer function of the given model. We follow the approach presented in~\cite{Schwerdtner2020} and adapted in~\cite{Schwerdtner2021Ident} to parameterize a pH system with a feedthrough term. The concept behind~\cite{Schwerdtner2020} is to construct the skew-symmetric and positive semi-definite parts of the  realization matrices of the ROM from strictly upper triangular matrices and upper triangular matrices, respectively. These triangular matrices are parameterized using the functions $\vtu(\cdot)$ (or $\vtsu(\cdot)$), which map a vector row-wise to an appropriately sized (strictly) upper triangular matrix. The function $\vtf_m(\cdot)$ is the standard reshape operation that maps a vector in $\R^{n \cdot m}$ to a matrix in $\R^{n \times m}$. These functions are explained in detail in~\cite[Definition~3.1]{Schwerdtner2020}.

\begin{lem}\label{lem:PHParam}
  Let $\theta \in \R^{n_\theta}$ be a parameter vector partitioned as $\theta = \begin{bmatrix}
    \theta_J^\T, \theta_W^\T, \theta_Q^\T, \theta_G^\T, \theta_N^\T
  \end{bmatrix}^\T$,
  with
  $\theta_J \in \R^{\dimx(\dimx-1)/2}$,
  $\theta_W \in \R^{(\dimx+\dimu)(\dimx+\dimu+1)/2}$,
  $\theta_Q \in \R^{\dimx(\dimx+1)/2}$,
  $\theta_G \in \R^{\dimx \cdot \dimu}$, and
  $\theta_N \in \R^{\dimu(\dimu-1)/2}$.
  Furthermore, define the matrix-valued functions
  \begin{subequations}
    \begin{align}
      \rJ(\theta) &:= \vtsu{(\theta_J)}^\T - \vtsu(\theta_J),\\
      \rW(\theta) &:= \vtu(\theta_W)\vtu{(\theta_W)}^\T  \label{eq:Wconstruction},\\
      \rR(\theta) &:= \begin{bmatrix} I_{\dimx} & 0 \end{bmatrix} \rW(\theta) \begin{bmatrix} I_{\dimx} & 0 \end{bmatrix}^\T,\\
      \rP(\theta) &:= \begin{bmatrix} I_{\dimx} & 0 \end{bmatrix} \rW(\theta) \begin{bmatrix} 0 & I_{\dimu} \end{bmatrix}^\T,\\
      \rS(\theta) &:= \begin{bmatrix} 0 & I_{\dimu} \end{bmatrix} \rW(\theta) \begin{bmatrix} 0 & I_{\dimu} \end{bmatrix}^\T, \\
       \rQ(\theta) &:= \vtu(\theta_Q)\vtu{(\theta_Q)}^\T,\\
      \rG(\theta) &:= \vtf_m(\theta_G),\\
      \rN(\theta) &:= \vtsu{(\theta_N)}^\T - \vtsu(\theta_N).
    \end{align}\label{eq:PHParamMatrices}
  \end{subequations}
  Then the parametric system
  \begin{align}
    \label{eq:pHParam}
    \pHsysr(\theta):
    \begin{cases}
      \!\begin{aligned}
        \dot{x}_r(t)= &(\rJ(\theta)-\rR(\theta))\rQ(\theta)\rx(t) \\ &\quad \quad + (\rG(\theta)-\rP(\theta)) u(t),
      \end{aligned}\\
      \!\begin{aligned}
        \ry(t)= &{(\rG(\theta)+\rP(\theta))}^\T \rQ(\theta) \rx(t)\\ &\quad \quad+ (\rS(\theta)-\rN(\theta)) u(t),
      \end{aligned}\\
    \end{cases}
  \end{align}
  is a pH-DAE (with $E_r=I_r$) as in Definition~\ref{def:pHsys}.
\end{lem}

The transfer function of a small-scale parametric system $\pHsysr(\theta)$ as in \eqref{eq:pHParam}, 
is denoted by $\pHtfr(\cdot, \theta)$.
The following comments motivate our choice for the parameterization.

\begin{enumerate}[(i)]
  \item Lemma~\ref{lem:PHParam} only allows the construction of pH-ODEs. Nevertheless, we can use the parameterization to approximate any pH-DAE with index one because we only aim to approximate the transfer function of the given system. Since transfer functions of descriptor systems with index one only have a constant polynomial part, we can approximate the effect of the algebraic constraints by tuning the feedthrough  terms $\rS(\theta)$ and $\rN(\theta)$ appropriately.
  \item Note that we can also rewrite the resulting pH-ODE as an implicit pH-DAE (with an $E_r$ term but without a $Q_r$ term), by a change of variables
  $\tilde x_r(t) := Q_r x_r(t)$ if $Q_r$ is positive definite.
  \item We ensure the positive semi-definiteness of the passivity matrix from Definition~\ref{def:pHsys}
    for all $\theta \in \R^{n_\theta}$ by ensuring that $\rW(\theta)$ is positive semi-definite for all $\theta \in \R^{n_\theta}$. The positive semi-definiteness of $\rQ(\theta)$ in combination with the skew-symmetry of $\rJ(\theta)$ also ensures $Q_r(\theta)^\T J_r(\theta) Q_r(\theta) = -Q_r(\theta)^\T J_r(\theta)^\T Q_r(\theta)$ for all $\theta \in \R^{n_\theta}$.
\end{enumerate}

In the following, we present two methods for tuning the parameter vector $\theta$ to obtain either an $\hinf$ or an $\htwo$ approximation for a given FOM\@.

\subsection{$\hinf$ approximation}%
\label{sub:ourmethod_hinf}

The algorithm presented in \cite{Schwerdtner2020} to obtain a good $\hinf$ approximation of a pH-FOM with transfer function $\pHtf$ using a parametrized low-order system with transfer function $\pHtfr( \cdot, \theta)$ is based on minimizing the objective function
\begin{align}
  \begin{split}
  \pushleft{\loss(\theta;\pHtf,\gamma,\mathcal{S}) :=} \\ \pushright{\frac{1}{\gamma}\sum\limits_{s_i\in \mathcal{S}}
  \sum\limits_{j=1}^{m}{\left({\left[\sigma_j \left(\pHtf(s_i)-\pHtfr(s_i,\theta)\right)-\gamma\right]}_+\right)}^2}
\end{split}
  \label{eq:loss}
\end{align}
with respect to $\theta$, where 
\begin{align*}
  {[ \cdot ]}_+:  \R \rightarrow [0,\infty), \quad x \mapsto 
  \begin{cases}
    x & \text{if } x\ge 0,\\
    0 & \text{if } x<0,
  \end{cases}
\end{align*}
for decreasing values of $\gamma > 0$. Here $\mathcal{S} \subset \ri \R$ is a set of sample points, at which the original and reduced transfer functions are evaluated and $\sigma_j(\cdot)$ denotes the $j$-th singular value of its matrix argument.

The justification for using 
$\loss$ 
as a surrogate for the $\hinf$ error is that $\loss(\cdot;H,\gamma,\mathcal{S})$ attains its global minimum (at zero), when all singular values of the error transfer function at all sample points $s_i$ are below the threshold~$\gamma$. Therefore, if the sample points are chosen appropriately, then $\loss(\theta;\pHtf,\gamma,\mathcal{S}) = 0$ is an indication for ${{\|\pHtf - \pHtf_r(\cdot,\theta)\|}_{\cH_\infty} \le \gamma}$. In \cite{Schwerdtner2021}, an adaptive sampling procedure is introduced which ensures that the sample points are appropriately distributed along the imaginary axis based on the given FOM and status of the optimization.

In this work, we use Algorithm~\ref{alg:bisection} to determine an $\hinf$ approximation via a bisection procedure, which determines the minimal value for $\gamma$ (up to a relative tolerance $\epsilon_1$) at which a minimization of $\loss$ with respect to $\theta$ terminates at zero. The tolerance $\epsilon_2$ that is used in line $6$ of the algorithm is the maximum value of $\loss$, which is still numerically interpreted as zero, such that $\gamma$ is reduced in the subsequent bisection step. The sample points are updated after each bisection step because the adaptive sampling update rule introduced in~\cite{Schwerdtner2021} depends on the current value of $\gamma$.

\begin{algorithm}[tbh]
  \LinesNumbered
  \SetAlgoLined
  \DontPrintSemicolon
  \SetKwInOut{Input}{Input}\SetKwInOut{Output}{Output}
  \Input{FOM transfer function $\pHtf\in\mathcal{RH}_{\infty}^{\enu \times \enu}$,
  initial ROM transfer function $\pHtfr( \cdot, \theta_0) \in \rhinf^{\enu \times \enu}$ with parameter $\theta_0 \in \R^{n_\theta}$,
  initial sample point set $\mathcal{S} \subset \mathrm{i}\R$,
  upper bound $\gamma_{\rm u} > 0$,
  bisection tolerance $\varepsilon_1 > 0$,
  termination tolerance $\varepsilon_2 > 0$
}
  \Output{Reduced pH-ODE of order $r$}
  Set $j:=0$ and $\gamma_{\rm l}:=0$.\;
  \While{$(\gamma_{\rm u}-\gamma_{\rm l})/(\gamma_{\rm u}+\gamma_{\rm l}) > \varepsilon_1$}{
    Set $\gamma:=(\gamma_{\rm u}+\gamma_{\rm l})/2$.\;
    Update the sample set $\mathcal{S}$ using~\cite[Alg.~3.1]{Schwerdtner2021}. \;
    Solve the minimization problem $\alpha := \min_{\theta\in \R^{n_\theta}}\loss(\theta;\pHtf,\gamma,\mathcal{S})$ with minimizer $\theta_{j+1} \in \R^{n_\theta}$, initialized at $\theta_j$. \;
  \eIf{$\alpha > \varepsilon_2$}{
    Set $\gamma_{\rm l}:=\gamma$.\;
    }{
    Set $\gamma_{\rm u}:=\gamma$.\;
  }
  Set $j:=j+1$.
  }
  Construct the ROM with $\theta_j$ as in Lemma~\ref{lem:PHParam}.
  \caption{SOBMOR-$\hinf$}
  \label{alg:bisection}
\end{algorithm}

The benefits of using this approach instead of directly minimizing the $\hinf$ norm were discussed in detail in~\cite[Remark~3.3]{Schwerdtner2020}. The main reasons for using~\eqref{eq:loss} instead of the $\hinf$ norm are the differentiability of~$\loss$ with respect to~$\theta$, the local convergence of the method, and the prohibitive computational costs as well as reliability issues of the $\hinf$ norm computation (for the large-scale error system) inside an optimization loop.

\subsection{$\htwo$ approximation}%
\label{sub:ourmethod_htwo}
In order to obtain a finite $\htwo$ error, the polynomial parts of the FOM and the ROM transfer function must be equal. When using projection-based methods on ODE models this \emph{feedthrough matching} is automatic. In the DAE case, this is typically obtained by preserving the algebraic part, i.\,e., by including the null-space of the $E$-matrix in the projection matrices. For systems with multiple algebraic constraints this is undesirable, and a reduction of the subsystem corresponding to the  algebraic constraints may be necessary, see, for instance, \cite{MehS05}.

Another remedy (used for pH-DAEs with index one in \cite{BeaGM21}) is to compute the polynomial part $\pHtf_{\rm p}(s) \equiv \mathcal{D}_0$ of the FOM before the reduction and include it in the feedthrough terms $\rS$ and $\rN$ of the ROM, which we will use here as well. The direct computation of $\mathcal{D}_0$ (see Section \ref{sec:preliminaries}) may, however, require transformations of the FOM. Alternatively, $\mathcal{D}_0$ can also be estimated by sampling the transfer function of the FOM at sufficiently large $s\in\C$ in an iterative manner as proposed in \cite{Schwerdtner2020a}. We then decompose it in its symmetric and skew-symmetric part, respectively, i.\,e., 
\begin{align}
  S_0 &:= \frac{1}{2}\left(\mathcal{D}_0^\T+\mathcal{D}_0\right), \label{eq:S0} \\ 
  N_0 &:= \frac{1}{2}\left(\mathcal{D}_0^\T-\mathcal{D}_0\right). \label{eq:N0} 
\end{align}
If the ROM is parameterized as in Lemma~\ref{lem:PHParam} then the $\cH_2$ error ${\Vert \pHtf - \pHtfr(\cdot,\theta)\Vert}_{\cH_2}$ is only well-defined if we have that 
\begin{align}\label{eq:SN_match}
	\rS(\theta)  &= S_0,\quad \rN(\theta) = N_0,
\end{align}
since otherwise $\pHtf-\pHtfr(\cdot,\theta) \notin \mathcal{RH}_{2}^{m \times m}$. Consequently,  we first have to fix all parameters in  $\theta_N$ such that ${\rN(\theta) = N_0}$. Since we indirectly parameterize $\rS(\theta)$ via $\theta_W$, we first analyze which parameters in $\theta_W$ have an impact on $\rS(\theta)$. For this, consider a separation of ${\theta_W \in \R^{n_W}}$ into ${\theta_W =: \mat{\theta_{W_1}^{\T}, \theta_{W_2}^{\T}}^\T}$, where ${\theta_{W_1} \in \R^{n_W - m(m+1)/2}}$ and ${\theta_{W_2} \in \R^{m(m+1)/2}}$. Then we can decompose 
\begin{align}
	\vtu(\theta_W) = \mat{\Xi_1 & \Xi_2 \\ 0 & \Xi_3},
\end{align}
where the matrices $\Xi_1, \,\Xi_2$ depend only on $\theta_{W_1}$, and $\Xi_3$ depends only on~$\theta_{W_2}$. Consequently, ${\rS(\theta) = \Xi_3\Xi_3^\T}$ only depends on~$\theta_{W_2}$ and we can set $\theta_{W_2}$ such that ${\rS(\theta) = S_0}$. The remaining parameters $\theta_{W_1}$ may still be subject to optimization and it holds that ${W(\theta) \geq 0}$ for all ${\theta_{W_1} \in \R^{n_W - m(m+1)/2}}$. Consequently, for minimizing the $\cH_2$ error, the parameter vector which is subject to optimization reduces to $\theta := \begin{bmatrix}
	\theta_J^\T, \theta_{W_1}^\T, \theta_Q^\T, \theta_{\pG}^\T
\end{bmatrix}^\T.$
\begin{rem}
  Note that $S_0 \geq 0$ always holds since the implicit pH-ODE~\eqref{eq:FOM_decpl} has the same transfer function $\pHtf$ as the original pH-DAE and is therefore positive real, i.\,e., we have that
  \begin{equation*}
    \Phi(\mathsf{i}\omega) := \pHtf(\mathsf{i}\omega) + \pHtf(\mathsf{i}\omega)^{\mathsf{H}} \geq 0, 
  \end{equation*}
  for all $\omega \in \R$ and consequently, $\lim_{\omega\rightarrow\infty}\Phi(\mathsf{i}\omega) = 2S_0\geq 0$.
\end{rem} 

Now we can formulate the $\cH_2$ optimization problem in the pole-residue framework originally proposed in \cite{Beattie2009a} for unstructured LTI systems and extended to pH-ODE systems in \cite{Moser2020}. Assume that
$(\rJ(\theta)-\rR(\theta))\rQ(\theta)$ is diagonalizable and consider the spectral decomposition
\begin{equation} \label{eq:red_evp}
	(\rJ(\theta)-\rR(\theta))\rQ(\theta) Z(\theta) = Z(\theta)\Lambda(\theta),
\end{equation}
where $\Lambda(\theta) = \text{diag}(\reig_1(\theta),\,...\,,\reig_r(\theta))$ and $Z(\theta)$ contains the right eigenvectors as columns. If the eigenvalues $\reig_i(\theta)~\in~\C$, $i=1,\,\ldots,\,r$ are simple, then the transfer function $\pHtfr(\cdot,\theta)$ may be represented by the partial fraction expansion ${\pHtfr(s,\theta) = \sum_{i=1}^{r} \frac{\rl_i(\theta)\rr_i(\theta)^\T}{s-\reig_i(\theta)}+\rS(\theta)-\rN(\theta)}$, where ${\rl_i(\theta),\,\rr_i(\theta)\in \C^m}$ with
\begin{align*}
	\rl_i(\theta) &= (\rG(\theta)+\rP(\theta))^\T\rQ(\theta) Z(\theta) e_i, \\
	\rr_i(\theta) &= (\rG(\theta)-\rP(\theta))^\T Z(\theta)^{-\T}e_i,
\end{align*}
and where $e_i$ denotes the $i$-th standard basis vector of $\R^r$. Assuming that $\eqref{eq:SN_match}$ holds, we have that	
\begin{equation}
	\begin{aligned}
		{\Vert \pHtf-\pHtfr(\cdot,\theta) \Vert}_{\cH_2}^2 = &{\Vert \pHtf_{\rm sp} \Vert}_{\cH_2}^2 \\	
				& - 2\sum_{i=1}^{r}\rl_i(\theta)^\T \pHtf_{\rm sp}(-\reig_i(\theta))\rr_i(\theta) \\ 
					& + \sum_{j,k=1}^{r} \frac{\rl_j(\theta)^\T\rl_k(\theta)\rr_k(\theta)^\T\rr_j(\theta)}{-\reig_j(\theta)-\reig_k(\theta)},
	\end{aligned}
\end{equation}
as shown in \cite[Theorem 2.1]{Beattie2009a}. Since ${\Vert \pHtf_{\rm sp} \Vert}_{\cH_2}^2$ does not depend on $\theta$, it can be neglected in the optimization. Consequently, we define the objective functional
\begin{align*}
\losst(\theta;\pHtf) &:={\Vert \pHtf-\pHtfr(\cdot,\theta) \Vert}_{\cH_2}^2 - {\Vert \pHtf_{\rm sp} \Vert}_{\cH_2}^2 \\
&:= \hat{\mathcal{F}}\Big(\big[\rl_1(\theta)^\T,\ldots ,\rl_r(\theta)^\T,\rr_1(\theta)^\T,\ldots\\ & \quad\quad\quad\quad \quad\rr_r(\theta)^\T,\reig_1(\theta)^\T,\ldots,\reig_r(\theta)^\T\big]^\T\Big) \\
&= (\hat{\mathcal{F}} \circ q)(\theta),
\end{align*}
where
\begin{align*}
	q(\theta) &:= [\rl_1(\theta)^\T,\ldots,\rl_r(\theta)^\T,\rr_1(\theta)^\T,\ldots \\ & \quad\quad\quad\quad\rr_r(\theta)^\T,\reig_1(\theta),\ldots,\reig_r(\theta)]^\T \in \C^{n_q}.
\end{align*}
This functional can be evaluated efficiently because it only requires the solution of the reduced-order eigenvalue problem in \eqref{eq:red_evp} as well as $r$ evaluations of $\pHtf_{\rm sp}$ at $-\reig_i(\theta)$. 
The eigenvalues $\reig_i(\theta)$ and rank-one residues $\rl_i(\theta)\rr_i(\theta)^\T$ are functions of the parameter vector $\theta$.
If $\bar{\theta} \in \R^{n_\theta}$ is chosen such that all eigenvalues are simple, then $\mathcal{F}$ is differentiable in a neighborhood of $\bar{\theta}$. 
Its derivative is obtained by applying the chain rule, i.\,e., with the differentiation operator $\text{D}$ we obtain
\begin{equation*}
 \text{D}\losst(\bar{\theta}) = \left(\nabla \losst(\bar{\theta})\right)^\T = \text{D}\hat{\mathcal{F}}(q(\bar{\theta}))\cdot \text{D}q(\bar{\theta}),
\end{equation*}
with
\begin{multline*}
	\text{D}\hat{\losst}(q(\bar{\theta})) = \left[\text{D}_{b_1}\hat{\losst}(q(\bar{\theta})),\ldots,\text{D}_{b_r}\hat{\losst}(q(\bar{\theta})), \ldots \right. \\ \left. \text{D}_{c_1}\hat{\losst}(q(\bar{\theta})),\ldots, \text{D}_{c_r}\hat{\losst}(q(\bar{\theta})), \ldots\right. \\ \left.\text{D}_{\lambda_1}\hat{\losst}(q(\bar{\theta})),\ldots,\text{D}_{\lambda_r}\hat{\losst}(q(\bar{\theta}))\right]\in\C^{1\times n_q},
\end{multline*}
and 
\begin{align*}
	\text{D}q(\bar{\theta}) &= \left[\text{D}_{\theta_1}{q(\bar{\theta})},\ldots,\text{D}_{\theta_{n_\theta}}{q(\bar{\theta})}\right]\in\C^{n_q\times n_{\theta}}.
\end{align*}
For all $i=1,\,\ldots,\,r$ it holds that
\begin{align*}
    \text{D}_{\rr_i} \hat{\losst}(q(\bar{\theta})) &= 2\rl_i(\bar{\theta})^\T\big(\pHtfr(-\reig_i(\bar{\theta}))-\pHtf(-\reig_i(\bar{\theta}))\big), \\
    \text{D}_{\rl_i} \hat{\losst}(q(\bar{\theta})) &= 2\rr_i(\bar{\theta})^\T\big(\pHtfr(-\reig_i(\bar{\theta}))-\pHtf(-\reig_i(\bar{\theta}))\big)^\T, \\
	\text{D}_{\lambda_i} \hat{\losst}(q(\bar{\theta})) &= -2\rl_i(\bar{\theta})^\T\big(\pHtfr^{\prime}(-\reig_i(\bar{\theta}))-\pHtf^{\prime}(-\reig_i(\bar{\theta}))\big)\rr_i(\bar{\theta}),
\end{align*}
and we refer  to~\cite{Moser2020} for the differentiation of $\text{D} q$. 
\begin{rem}
	The partial derivatives in $\mathrm{D}q(\bar{\theta})$ may be computed efficiently with block-wise expressions. For instance, the derivative $\mathrm{D}_{\theta_G}\rl_i(\bar{\theta}) \in \C^{m \times r\cdot m}$ can be computed as
	\begin{align*}
		\text{D}_{\theta_G}\rl_i(\bar{\theta}) &=  \mat{z_i(\bar{\theta})^\T\rQ(\bar{\theta}) & 0 & 0 \\ 0 & \ddots & 0 \\ 0 & 0 & z_i(\bar{\theta})^\T\rQ(\bar{\theta})} \\ &= I_m \otimes  z_i(\bar{\theta})^\T\rQ(\bar{\theta}),
	\end{align*}
	where $z_i(\bar{\theta})\in\C^r$ denotes the $i$-th column in $Z(\bar{\theta})$. This is also the case for more complex derivatives that involve the differentiation of the eigenvalue problem in \eqref{eq:red_evp}.
\end{rem}
Here, we highlight some important advantages of the pole-residue framework compared to recently proposed methods that are formulated in the Lyapunov framework (see~\cite{Sato2018,Jiang2019}), in particular for pH-DAEs. These methods require the solution of large-scale Lyapunov equations for the evaluation of $\losst$ and its gradient. Currently no structure-preserving Lyapunov-based methods exist for pH-DAEs; see \cite{SchM22_ppt} for new Lyapunov-based formulations of pH-DAEs. If the strictly proper part of the transfer function can be easily decoupled from the constant polynomial part, then for pH-DAEs with index one, the existing methods for pH-ODEs may be applied to this part.
However, if the splitting into the strictly proper and polynomial part has first to be computed via a factorization method, then the sparsity patterns of the original pH-DAE may be lost
which complicates the repetitive solution of Lyapunov equations for these systems in the large-scale setting. We highlight that the pole-residue framework only requires evaluations of $\pHtf_{\rm sp}$. Since we have that
\begin{equation}
	\pHtf_{\rm sp}(s) = \pHtf(s) - (S_0-N_0),
\end{equation}
for all $s\in\C$, we may work directly with the sparse matrices of the original pH-DAE and do not require the solution of large-scale Lyapunov equations.
\begin{algorithm}[tbh]
	\LinesNumbered
	\SetAlgoLined
	\DontPrintSemicolon
	\SetKwInOut{Input}{Input}\SetKwInOut{Output}{Output}
	\Input{FOM transfer function $\pHtf\in\mathcal{RH}_{\infty}^{\enu \times \enu}$, reduced order $r\in\mathbb{N}$.
	}
	\Output{Reduced pH-ODE of order $r$}
	Compute $S_0,\,N_0$ as in \eqref{eq:S0}--\eqref{eq:N0}.\;
	Initialize $\theta_0$ s.t. $\rS(\theta_0)=S_0$, $\rN(\theta_0)=N_0$.\;
	Solve \vspace{-12pt}\begin{align*}
	    \theta_{\rm fin} &= \argmin\limits_{\substack{\theta\in \R^{n_\theta}}} \losst(\theta; \pHtf) \\
	    &\text{s.t. } \rS(\theta)=S_0,\, \rN(\theta)=N_0.
	\end{align*}\vspace{-18pt} \; 
	Construct the ROM with $\theta_{\rm fin}$ as in Lemma \ref{lem:PHParam}.\;
	\caption{PROPT-$\mathcal{H}_2$}
	\label{alg:H2OPT}
\end{algorithm}

Since the $\cH_2$ optimization problem is non-convex, the choice of the initial parameter vector $\theta_0$ will generally impact the fidelity of the final ROM obtained by Algorithm~\ref{alg:H2OPT}. Simple initialization strategies are, for instance, choosing $\theta_0$ randomly or using IRKA-PH (see~\cite{Moser2020, Sato2018}), which generally converges very quickly. Here, we propose another approach that may use \emph{unstructured} ROMs for initialization which is based on the following parameterization.
\begin{lem}\label{lem:pHInit}
	 Let $(\widetilde{A},\widetilde{B},\widetilde{C},\widetilde{D})$ be a ROM of state-space dimension $r$ such that $\widetilde{D} = S_0-N_0$ and such that $\widetilde{A}$ has all its eigenvalues in the open left half of the complex plane. Let  $\theta_G \in\R^{r\cdot m}$ and $\theta_K \in\R^{r\cdot p}$ be two parameter vectors and define the matrix-valued functions
	 \begin{align*}
	    \iG(\theta_G) &:= \vtf_m(\theta_G), \\
	 	\rK(\theta_K) &:= \vtf_q(\theta_K).
	 \end{align*}
 	 Let $\iQ(\theta_K)>0$ solve the Lyapunov equation
 	 \begin{equation}
 	 	\widetilde{A}^\T\iQ(\theta_K)+\iQ(\theta_K)\widetilde{A}+\rK(\theta_K)\rK(\theta_K)^\T = 0,
 	 \end{equation}
  	and define 
  	\begin{align*}
  		\iJ(\theta_K)&=\frac{1}{2}\left(\widetilde{A} \iQ(\theta_K)^{-1}-\iQ(\theta_K)^{-1}\widetilde{A}^\T\right), \\
  		\iR(\theta_K)&=-\frac{1}{2}\left(\widetilde{A} \iQ(\theta_K)^{-1}+\iQ(\theta_K)^{-1}\widetilde{A}^\T\right).
  	\end{align*}
  	Then the parametric system
  	\begin{align}
  		\label{eq:pHParamInit}
  		\pHsysr(\theta_G,\theta_K):
  		\begin{cases}
  			\!\begin{aligned}
  				\dot x_r(t)= &\big(\iJ(\theta_K)-\iR(\theta_K))\iQ(\theta_K\big)\rx(t) \\ &\quad \quad + \iG(\theta_G) u(t),
  			\end{aligned}\\
  			\!\begin{aligned}
  		 	\ry(t)= &\iG(\theta_G)^\T \iQ(\theta_K) \rx(t)\\ &\quad \quad+ (S_0-N_0) u(t)
  			\end{aligned}\\
  		\end{cases}
  	\end{align}
  	is a pH-ODE system with
  	\begin{equation*}
  		\big(\iJ(\theta_K)-\iR(\theta_K)\big)\iQ(\theta_K) = \widetilde{A}.
  	\end{equation*}  	 
\end{lem}
Let $\widetilde{H}$ denote the transfer function of the (possibly unstructured) ROM $\big(\widetilde{A},\widetilde{B},\widetilde{C},\widetilde{D}\big)$ with ${\widetilde{\pHtf}(s)=\sum\limits_{i=1}^{r} \frac{\widetilde{\rl}_i\widetilde{\rr}_i^\T}{s-\widetilde{\reig}_i}+S_0-N_0}$ and ${\widetilde{\rl}_i,\,\widetilde{\rr}_i\in \C^m}$. 

Based on the parameterization in Lemma \ref{lem:pHInit}, we can then compute an initial pH model by minimizing the weighted sum of squared errors between the residuals in the Frobenius norm,  i.\,e.,  
\begin{equation*}
	\losst_0(\theta_G,\theta_K) := \sum_{i=1}^{r} \frac{1}{|\widetilde{\reig}_i|} \left\Vert \widetilde{\rl}_i\widetilde{\rr}_i^\T - \rl_i(\theta_G,\theta_K)\rr_i(\theta_G,\theta_K)^\T\right\Vert_{\rm F}^2,
\end{equation*}
where
\begin{align*}
	\rl_i(\theta_G,\theta_K) &= \iG(\theta_G)^\T\iQ(\theta_K) \widetilde{Z} e_i, \\
	\rr_i(\theta_G,\theta_K) &= \iG(\theta_G)^\T \widetilde{Z}^{-\T}e_i,
\end{align*}
for $i=1,\,\ldots,\,r$ and $\widetilde{Z}$ is, 
again under a diagonalizability assumption, obtained from the spectral decomposition
 \begin{equation*} 
	\widetilde{A} \widetilde{Z} = \widetilde{Z}\widetilde{\Lambda},
\end{equation*}
with $\widetilde{\Lambda} = \text{diag}(\widetilde{\reig}_1,\,...\,,\widetilde{\reig}_r)$.

Note that the computation of the gradient of $\losst_0$ is very simple, since it does not involve a differentiation of the eigenvalues or eigenvectors. While the partial gradients of $\rl_i(\cdot)$ and $\rr_i(\cdot)$ with respect to $\theta_G$ are straightforward, the partial gradients of $\iQ(\cdot)$ with respect to the $l$-th entry in $\theta_K$ is the solution of the (reduced-order) Lyapunov equation
\begin{multline*}
	\widetilde{A}^\T \frac{\partial\iQ(\theta_K)}{\partial \theta_{K,l}} + \frac{\partial\iQ(\theta_K)}{\partial \theta_{K,l}}\widetilde{A} \\ + \vtf_q(e_l)\rK(\theta_K)^\T+\rK(\theta_K)\vtf_q(e_l)^\T= 0,
\end{multline*}
where $e_l$ denotes the $l$-th standard basis vector of $\R^{r\cdot p}$. As the number of optimization parameters is reduced to $r(p+m)$, this initialization generally converges very quickly. In combination with Algorithm~\ref{alg:H2OPT}, this enables a \emph{two-step} approach with a more restrictive (yet simpler) pre-optimization of only the residuals and a subsequent (more complex) optimization of all system matrices.
\begin{rem}
	Note that the sample-based SOBMOR method can be tuned to compute a ROM with small $\htwo$ error as well. Instead of using $\loss$ in conjunction with the bisection method outlined in Algorithm~\ref{alg:bisection}, the integral in~\eqref{eq:htwonorm} can be approximated by means of an adaptive quadrature rule (see \cite[Algorithm~1]{Gonnet2010} for a template method). In this way, it is possible to compute the $\htwo$ error and its gradient with respect to the free ROM parameters in terms of the error transfer function at specific sample points. This makes it possible to use the same optimization techniques as in~\cite{Schwerdtner2020}. In particular, \cite[Theorem~3.1]{Schwerdtner2020} for the gradient computation can be reused. 
	
	We denote this method by SOBMOR-$\htwo$. Further details regarding the implementation of the adaptive integration are provided in the Appendix.
\end{rem}

\section{Numerical examples}%
\label{sec:numerical_examples}

\begin{figure}[tpb]
	\centering
	\includegraphics[width=0.45\textwidth]{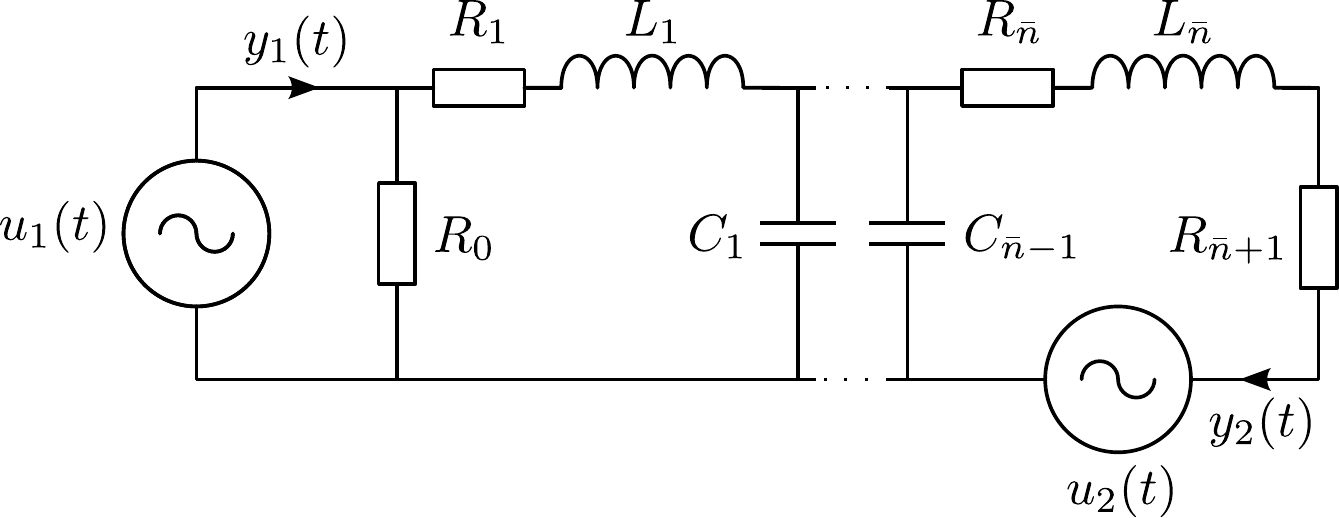}
	\caption{RCL ladder network with two voltage sources.}\label{fig:RCL_Ladder}
\end{figure}
To evaluate our approach, we consider different variants of an RCL ladder network as in Figure~\ref{fig:RCL_Ladder}. Ladder networks are often used as surrogate models for transmission lines in power networks, see \cite{GreBR13}. The number of loops in the network is denoted by $\bar{n}$. In the following, we consider two configurations of the system: a multiple-input multiple-output (MIMO) version, where the inputs are the voltages of both voltage sources and the outputs are the currents as shown in Figure \ref{fig:RCL_Ladder}. In the single-input single-output (SISO) configuration, we replace the second voltage source by a wire and only consider the input-to-output behaviour from $u_1(\cdot)$ to $y_1(\cdot)$.

\begin{figure}[t]
  \centering
  \input{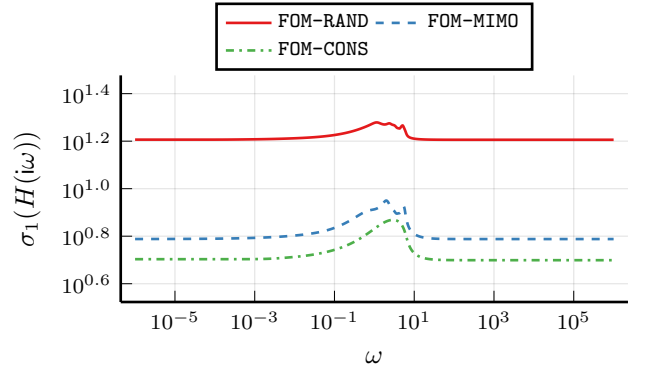}
  \caption{Maximal singular values of the FOM transfer functions.}\label{fig:FOMs}
\end{figure}

Modeling of RCL circuits as depicted in Figure~\ref{fig:RCL_Ladder} via the lumped-element approach described in~\cite{Freund2011} directly leads to pH-DAE models with index one. For further details about the model we refer to the software package \textsf{PortHamiltonianBenchmarkSystems}\footnote{\url{https://algopaul.github.io/PortHamiltonianBenchmarkSystems/RclCircuits/}}, which we use to generate three different RCL circuits.
Our first model, \texttt{FOM-CONS}, contains 100 loops and the inductances, resistances, and capacities are the same in each loop. For \texttt{FOM-RAND} and \texttt{FOM-MIMO}, the resistances are chosen randomly to obtain a more complex model, that we expect to be harder to reduce. The models contain 500 and $10\,000$ loops, respectively.
The key dimensions of all considered FOMs are sumarized in Table~\ref{tab:systems}. In Figure~\ref{fig:FOMs}, the transfer functions of all three FOMs are displayed\footnote{All FOM system matrices are available at \url{https://doi.org/10.5281/zenodo.6497076}}. As expected, \texttt{FOM-CONS} results in just one smooth peak, while \texttt{FOM-RAND} and \texttt{FOM-MIMO} have several and sharper peaks.

\begin{table}[t]
  \centering
  \caption{Dimensions of benchmark systems. The number of input/output pairs is denoted by $m$.}\label{tab:systems}
    \begin{tabular}{c|cccc}
         model name & $\bar n$ & $n$      & $m$ & $\rank E$\\ \hline
      \texttt{FOM-CONS}    & $100$    & $302$    & 1   & $199$ \\
      \texttt{FOM-RAND}    & $500$    & $1502$  & 1   & $999$\\
      \texttt{FOM-MIMO}    & $10\,000$ & $30\,004$ & 2   & $19\,999$\\
    \end{tabular}
\end{table}

We first report the $\hinf$ accuracy of the considered methods in Figures~\ref{fig:FOMRANDORIGComp} (a) and (c) for the systems \texttt{FOM-CONS} and \texttt{FOM-RAND}. It can be seen that SOBMOR-$\hinf$ achieves the highest $\hinf$ accuracy across both models and all reduced model orders. The second best overall accuracy is obtained by PRBT and XminBT, which invariably have similar $\hinf$ errors. The $\htwo$ methods PROPT-$\htwo$, XminIRKA, and IRKA-PH have the worst $\hinf$ performance, as it is to be expected. However, we note that (especially in Figure~\ref{fig:FOMRANDORIGComp} (a) there is a huge difference in terms of accuracy, when comparing IRKA-PH to all the other methods. In general the spread of accuracies is higher for the more complicated \texttt{FOM-RAND}\@. Here, the~$\hinf$ accuracy of PRBT and XminBT is sometimes even worse than that of PROPT-$\htwo$\@.

\begin{figure*}[htpb]
  \centering
    \begin{tabular}{cc}
  \begin{tikzpicture}[/tikz/background rectangle/.style={fill={rgb,1:red,1.0;green,1.0;blue,1.0}, draw opacity={1.0}}, show background rectangle]
\begin{axis}[point meta max={nan}, point meta min={nan}, legend cell align={left}, legend columns={3}, title={}, title style={at={{(0.5,1)}}, anchor={south}, font={{\fontsize{14 pt}{18.2 pt}\selectfont}}, color={rgb,1:red,0.0;green,0.0;blue,0.0}, draw opacity={1.0}, rotate={0.0}}, legend style={color={rgb,1:red,0.0;green,0.0;blue,0.0}, draw opacity={1.0}, line width={1}, solid, fill={rgb,1:red,1.0;green,1.0;blue,1.0}, fill opacity={1.0}, text opacity={1.0}, font={{\fontsize{8 pt}{10.4 pt}\selectfont}}, text={rgb,1:red,0.0;green,0.0;blue,0.0}, cells={anchor={center}}, at={(0.5, 1.05)}, anchor={south}}, axis background/.style={fill={rgb,1:red,1.0;green,1.0;blue,1.0}, opacity={1.0}}, anchor={north west}, xshift={1.0mm}, yshift={-1.0mm}, clip mode={individual}, width={0.45\textwidth}, height={0.25\textwidth}, scaled x ticks={false}, xlabel={\scriptsize{Reduced model order $r$}}, x tick style={color={rgb,1:red,0.0;green,0.0;blue,0.0}, opacity={1.0}}, x tick label style={color={rgb,1:red,0.0;green,0.0;blue,0.0}, opacity={1.0}, rotate={0}}, xlabel style={}, xmajorgrids={true}, xmin={1.46}, xmax={20.54}, xtick={{2.0,6.0,10.0,14.0,18.0}}, xticklabels={{$2$,$6$,$10$,$14$,$18$}}, xtick align={inside}, xticklabel style={font={{\fontsize{8 pt}{10.4 pt}\selectfont}}, color={rgb,1:red,0.0;green,0.0;blue,0.0}, draw opacity={1.0}, rotate={0.0}}, x grid style={color={rgb,1:red,0.0;green,0.0;blue,0.0}, draw opacity={0.1}, line width={0.5}, solid}, axis x line*={left}, x axis line style={color={rgb,1:red,0.0;green,0.0;blue,0.0}, draw opacity={1.0}, line width={1}, solid}, scaled y ticks={false}, ylabel={\scriptsize{$\mathcal{H}_\infty$ error}}, y tick style={color={rgb,1:red,0.0;green,0.0;blue,0.0}, opacity={1.0}}, y tick label style={color={rgb,1:red,0.0;green,0.0;blue,0.0}, opacity={1.0}, rotate={0}}, ylabel style={}, ymode={log}, log basis y={10}, ymajorgrids={true}, ymin={6.016335111704473e-8}, ymax={2.110293638084952}, ytick={{1.0,0.01,0.0001,1.0e-6}}, yticklabels={{$10^{0}$,$10^{-2}$,$10^{-4}$,$10^{-6}$}}, ytick align={inside}, yticklabel style={font={{\fontsize{8 pt}{10.4 pt}\selectfont}}, color={rgb,1:red,0.0;green,0.0;blue,0.0}, draw opacity={1.0}, rotate={0.0}}, y grid style={color={rgb,1:red,0.0;green,0.0;blue,0.0}, draw opacity={0.1}, line width={0.5}, solid}, axis y line*={left}, y axis line style={color={rgb,1:red,0.0;green,0.0;blue,0.0}, draw opacity={1.0}, line width={1}, solid}, colorbar={false}]
    \addplot[color={rgb,1:red,0.5961;green,0.3059;blue,0.6392}, name path={d9d9839e-d6dc-4a62-abb7-32e69d19429a}, draw opacity={1.0}, line width={1}, solid, mark={star}, mark size={3.0 pt}, mark repeat={1}, mark options={color={rgb,1:red,0.0;green,0.0;blue,0.0}, draw opacity={1.0}, fill={rgb,1:red,0.5961;green,0.3059;blue,0.6392}, fill opacity={1.0}, line width={0.75}, rotate={0}, solid}]
        table[row sep={\\}]
        {
            \\
            2.0  0.8398640453341348  \\
            4.0  0.09996025983546461  \\
            6.0  0.05077950335584149  \\
            8.0  0.0044329393294906825  \\
            10.0  0.0021293366179922137  \\
            12.0  0.0001440616144953424  \\
            14.0  5.468451083073914e-5  \\
            16.0  5.364459477051393e-6  \\
            18.0  1.6376631488392914e-6  \\
            20.0  1.6250083252075111e-6  \\
        }
        ;
    \addlegendentry {XminBT$ $}
    \addplot[color={rgb,1:red,0.8941;green,0.102;blue,0.1098}, name path={2ac466ba-9356-43da-9d97-158ec8ef3c77}, draw opacity={1.0}, line width={1}, solid, mark={square*}, mark size={3.0 pt}, mark repeat={1}, mark options={color={rgb,1:red,0.0;green,0.0;blue,0.0}, draw opacity={1.0}, fill={rgb,1:red,0.8941;green,0.102;blue,0.1098}, fill opacity={1.0}, line width={0.75}, rotate={0}, solid}]
        table[row sep={\\}]
        {
            \\
            2.0  1.2790550707326978  \\
            4.0  0.46073627112202353  \\
            6.0  0.15280179996631232  \\
            8.0  0.06265335464104474  \\
            10.0  0.006754958288584514  \\
            12.0  0.00021249929733045027  \\
            14.0  0.00011356210035486439  \\
            16.0  5.3024328726007844e-6  \\
            18.0  2.5854527500759174e-6  \\
            20.0  1.7803459880859873e-6  \\
        }
        ;
    \addlegendentry {XminIRKA$ $}
    \addplot[color={rgb,1:red,0.651;green,0.3373;blue,0.1569}, name path={ab0acbce-6761-4778-a770-92c020841054}, draw opacity={1.0}, line width={1}, solid, mark={triangle*}, mark size={3.0 pt}, mark repeat={1}, mark options={color={rgb,1:red,0.0;green,0.0;blue,0.0}, draw opacity={1.0}, fill={rgb,1:red,0.651;green,0.3373;blue,0.1569}, fill opacity={1.0}, line width={0.75}, rotate={0}, solid}]
        table[row sep={\\}]
        {
            \\
            2.0  0.8123497251462448  \\
            4.0  0.10484079984099783  \\
            6.0  0.04758085564803856  \\
            8.0  0.0041575377898128004  \\
            10.0  0.001969479580802992  \\
            12.0  0.00015162158333626913  \\
            14.0  5.051855984411604e-5  \\
            16.0  5.632665190809335e-6  \\
            18.0  1.6140242306516173e-6  \\
            20.0  1.5329378518874412e-6  \\
        }
        ;
    \addlegendentry {PRBT$ $}
    \addplot[color={rgb,1:red,1.0;green,0.498;blue,0.0}, name path={69302f1b-9f3a-4cb5-b1b1-137c6f9edeb2}, draw opacity={1.0}, line width={1}, solid, mark={triangle*}, mark size={3.0 pt}, mark repeat={1}, mark options={color={rgb,1:red,0.0;green,0.0;blue,0.0}, draw opacity={1.0}, fill={rgb,1:red,1.0;green,0.498;blue,0.0}, fill opacity={1.0}, line width={0.75}, rotate={180}, solid}]
        table[row sep={\\}]
        {
            \\
            2.0  1.0671362504032225  \\
            4.0  0.5337771831462381  \\
            6.0  0.31180490319371723  \\
            8.0  0.180807555556548  \\
            10.0  0.09353485989267224  \\
            12.0  0.051333002741814376  \\
            14.0  0.028429112537013305  \\
            16.0  0.015416371042563076  \\
            18.0  0.008106592790858411  \\
            20.0  0.003931676655057896  \\
        }
        ;
    \addlegendentry {IRKA-PH$ $}
    \addplot[color={rgb,1:red,0.302;green,0.6863;blue,0.2902}, name path={3810ba5a-a5c6-46af-bed7-a5ceaa8b6546}, draw opacity={1.0}, line width={1}, solid, mark={*}, mark size={3.0 pt}, mark repeat={1}, mark options={color={rgb,1:red,0.0;green,0.0;blue,0.0}, draw opacity={1.0}, fill={rgb,1:red,0.302;green,0.6863;blue,0.2902}, fill opacity={1.0}, line width={0.75}, rotate={0}, solid}]
        table[row sep={\\}]
        {
            \\
            2.0  0.31911473400751106  \\
            4.0  0.04909386548777388  \\
            6.0  0.016619052400327978  \\
            8.0  0.003076699882813701  \\
            10.0  0.0007706866197239086  \\
            12.0  9.931950957306188e-5  \\
            14.0  2.088092877211946e-5  \\
            16.0  2.5961725048299098e-6  \\
            18.0  4.933389280864454e-7  \\
            20.0  9.837167405253637e-8  \\
        }
        ;
    \addlegendentry {SOBMOR-$\mathcal{H}_\infty$}
    \addplot[color={rgb,1:red,0.2157;green,0.4941;blue,0.7216}, name path={2567e106-5855-4b51-bbfc-882a7acae275}, draw opacity={1.0}, line width={1}, solid, mark={diamond*}, mark size={3.0 pt}, mark repeat={1}, mark options={color={rgb,1:red,0.0;green,0.0;blue,0.0}, draw opacity={1.0}, fill={rgb,1:red,0.2157;green,0.4941;blue,0.7216}, fill opacity={1.0}, line width={0.75}, rotate={0}, solid}]
        table[row sep={\\}]
        {
            \\
            2.0  1.2906391837996507  \\
            4.0  0.10690882859742615  \\
            6.0  0.15789652341508872  \\
            8.0  0.01138041479748891  \\
            10.0  0.007510545323324757  \\
            12.0  0.00022804878446767765  \\
            14.0  0.00011699188668496981  \\
            16.0  5.220220689790122e-6  \\
            18.0  3.1547954289889556e-6  \\
            20.0  2.3280764387621984e-6  \\
        }
        ;
    \addlegendentry {PROPT-$\htwo$}
\end{axis}
\end{tikzpicture} & \begin{tikzpicture}[/tikz/background rectangle/.style={fill={rgb,1:red,1.0;green,1.0;blue,1.0}, draw opacity={1.0}}, show background rectangle]
\begin{axis}[point meta max={nan}, point meta min={nan}, legend cell align={left}, legend columns={3}, title={}, title style={at={{(0.5,1)}}, anchor={south}, font={{\fontsize{14 pt}{18.2 pt}\selectfont}}, color={rgb,1:red,0.0;green,0.0;blue,0.0}, draw opacity={1.0}, rotate={0.0}}, legend style={color={rgb,1:red,0.0;green,0.0;blue,0.0}, draw opacity={1.0}, line width={1}, solid, fill={rgb,1:red,1.0;green,1.0;blue,1.0}, fill opacity={1.0}, text opacity={1.0}, font={{\fontsize{8 pt}{10.4 pt}\selectfont}}, text={rgb,1:red,0.0;green,0.0;blue,0.0}, cells={anchor={center}}, at={(0.5, 1.05)}, anchor={south}}, axis background/.style={fill={rgb,1:red,1.0;green,1.0;blue,1.0}, opacity={1.0}}, anchor={north west}, xshift={1.0mm}, yshift={-1.0mm}, clip mode={individual}, width={0.45\textwidth}, height={0.25\textwidth}, scaled x ticks={false}, xlabel={\scriptsize{Reduced model order $r$}}, x tick style={color={rgb,1:red,0.0;green,0.0;blue,0.0}, opacity={1.0}}, x tick label style={color={rgb,1:red,0.0;green,0.0;blue,0.0}, opacity={1.0}, rotate={0}}, xlabel style={}, xmajorgrids={true}, xmin={1.46}, xmax={20.54}, xtick={{2.0,6.0,10.0,14.0,18.0}}, xticklabels={{$2$,$6$,$10$,$14$,$18$}}, xtick align={inside}, xticklabel style={font={{\fontsize{8 pt}{10.4 pt}\selectfont}}, color={rgb,1:red,0.0;green,0.0;blue,0.0}, draw opacity={1.0}, rotate={0.0}}, x grid style={color={rgb,1:red,0.0;green,0.0;blue,0.0}, draw opacity={0.1}, line width={0.5}, solid}, axis x line*={left}, x axis line style={color={rgb,1:red,0.0;green,0.0;blue,0.0}, draw opacity={1.0}, line width={1}, solid}, scaled y ticks={false}, ylabel={\scriptsize{$\mathcal{H}_2$ error}}, y tick style={color={rgb,1:red,0.0;green,0.0;blue,0.0}, opacity={1.0}}, y tick label style={color={rgb,1:red,0.0;green,0.0;blue,0.0}, opacity={1.0}, rotate={0}}, ylabel style={}, ymode={log}, log basis y={10}, ymajorgrids={true}, ymin={6.016335111704473e-8}, ymax={2.3691440906378167}, ytick={{1.0,0.01,0.0001,1.0e-6}}, yticklabels={{$10^{0}$,$10^{-2}$,$10^{-4}$,$10^{-6}$}}, ytick align={inside}, yticklabel style={font={{\fontsize{8 pt}{10.4 pt}\selectfont}}, color={rgb,1:red,0.0;green,0.0;blue,0.0}, draw opacity={1.0}, rotate={0.0}}, y grid style={color={rgb,1:red,0.0;green,0.0;blue,0.0}, draw opacity={0.1}, line width={0.5}, solid}, axis y line*={left}, y axis line style={color={rgb,1:red,0.0;green,0.0;blue,0.0}, draw opacity={1.0}, line width={1}, solid}, colorbar={false}]
    \addplot[color={rgb,1:red,0.5961;green,0.3059;blue,0.6392}, name path={4ce6e899-4e1b-49ab-9c92-bbb9c4ad0327}, draw opacity={1.0}, line width={1}, solid, mark={star}, mark size={3.0 pt}, mark repeat={1}, mark options={color={rgb,1:red,0.0;green,0.0;blue,0.0}, draw opacity={1.0}, fill={rgb,1:red,0.5961;green,0.3059;blue,0.6392}, fill opacity={1.0}, line width={0.75}, rotate={0}, solid}]
        table[row sep={\\}]
        {
            \\
            2.0  0.6538838735595841  \\
            4.0  0.1461081749242702  \\
            6.0  0.032063895115099356  \\
            8.0  0.005683955850135333  \\
            10.0  0.0014337429633624111  \\
            12.0  0.00021081376861001473  \\
            14.0  4.793620564427977e-5  \\
            16.0  7.755332399021324e-6  \\
            18.0  1.3161478579172374e-6  \\
            20.0  1.3047212850631054e-6  \\
        }
        ;
    \addlegendentry {XminBT$ $}
    \addplot[color={rgb,1:red,0.8941;green,0.102;blue,0.1098}, name path={b980f0cd-f033-453b-ab2d-f46a3e810bd8}, draw opacity={1.0}, line width={1}, solid, mark={square*}, mark size={3.0 pt}, mark repeat={1}, mark options={color={rgb,1:red,0.0;green,0.0;blue,0.0}, draw opacity={1.0}, fill={rgb,1:red,0.8941;green,0.102;blue,0.1098}, fill opacity={1.0}, line width={0.75}, rotate={0}, solid}]
        table[row sep={\\}]
        {
            \\
            2.0  0.53180268044148  \\
            4.0  0.09586967703072999  \\
            6.0  0.015657263680482734  \\
            8.0  0.0052511615879096965  \\
            10.0  0.0011163495713920198  \\
            12.0  0.00019540960619362862  \\
            14.0  4.0503688706619446e-5  \\
            16.0  7.2444144924638565e-6  \\
            18.0  9.171610406096189e-7  \\
            20.0  4.16126208868108e-7  \\
        }
        ;
    \addlegendentry {XminIRKA$ $}
    \addplot[color={rgb,1:red,0.651;green,0.3373;blue,0.1569}, name path={3b953184-d85f-497b-918b-5c24fa09a1ab}, draw opacity={1.0}, line width={1}, solid, mark={triangle*}, mark size={3.0 pt}, mark repeat={1}, mark options={color={rgb,1:red,0.0;green,0.0;blue,0.0}, draw opacity={1.0}, fill={rgb,1:red,0.651;green,0.3373;blue,0.1569}, fill opacity={1.0}, line width={0.75}, rotate={0}, solid}]
        table[row sep={\\}]
        {
            \\
            2.0  0.6858120724257601  \\
            4.0  0.15187534539230033  \\
            6.0  0.03423603967460752  \\
            8.0  0.005941260548763478  \\
            10.0  0.0015212531691354472  \\
            12.0  0.0002195070776342319  \\
            14.0  5.048716498521421e-5  \\
            16.0  8.061768019369066e-6  \\
            18.0  1.3787936858462582e-6  \\
            20.0  1.2628988992695961e-6  \\
        }
        ;
    \addlegendentry {PRBT$ $}
    \addplot[color={rgb,1:red,1.0;green,0.498;blue,0.0}, name path={c24a87a3-f18d-41d8-b0d9-9b2d2e9e14c9}, draw opacity={1.0}, line width={1}, solid, mark={triangle*}, mark size={3.0 pt}, mark repeat={1}, mark options={color={rgb,1:red,0.0;green,0.0;blue,0.0}, draw opacity={1.0}, fill={rgb,1:red,1.0;green,0.498;blue,0.0}, fill opacity={1.0}, line width={0.75}, rotate={180}, solid}]
        table[row sep={\\}]
        {
            \\
            2.0  1.5060286382358836  \\
            4.0  0.7080989901786512  \\
            6.0  0.40429182172274836  \\
            8.0  0.22833852574590646  \\
            10.0  0.11407893953689882  \\
            12.0  0.061860921522758856  \\
            14.0  0.03382304516499245  \\
            16.0  0.018059607732053572  \\
            18.0  0.009411942450718264  \\
            20.0  0.004520734650664615  \\
        }
        ;
    \addlegendentry {IRKA-PH$ $}
    \addplot[color={rgb,1:red,0.302;green,0.6863;blue,0.2902}, name path={3bf4693f-123b-44c6-9218-98d774b7a557}, draw opacity={1.0}, line width={1}, solid, mark={*}, mark size={3.0 pt}, mark repeat={1}, mark options={color={rgb,1:red,0.0;green,0.0;blue,0.0}, draw opacity={1.0}, fill={rgb,1:red,0.302;green,0.6863;blue,0.2902}, fill opacity={1.0}, line width={0.75}, rotate={0}, solid}]
        table[row sep={\\}]
        {
            \\
            2.0  0.5303912128947575  \\
            4.0  0.09556035301300327  \\
            6.0  0.015612207134097913  \\
            8.0  0.005231488339439383  \\
            10.0  0.0011115370085658307  \\
            12.0  0.00033186000603871506  \\
            14.0  4.034955601834643e-5  \\
            16.0  9.418957255372185e-6  \\
            18.0  4.626577732836335e-6  \\
            20.0  2.2101985080841075e-6  \\
        }
        ;
    \addlegendentry {SOBMOR-$\mathcal{H}_2$}
    \addplot[color={rgb,1:red,0.2157;green,0.4941;blue,0.7216}, name path={9c5ac194-4fa4-46ab-abec-3d82a48e7d36}, draw opacity={1.0}, line width={1}, solid, mark={diamond*}, mark size={3.0 pt}, mark repeat={1}, mark options={color={rgb,1:red,0.0;green,0.0;blue,0.0}, draw opacity={1.0}, fill={rgb,1:red,0.2157;green,0.4941;blue,0.7216}, fill opacity={1.0}, line width={0.75}, rotate={0}, solid}]
        table[row sep={\\}]
        {
            \\
            2.0  0.5303899877028492  \\
            4.0  0.1356087406019664  \\
            6.0  0.015612186586459261  \\
            8.0  0.005160543126538926  \\
            10.0  0.0011115368941180049  \\
            12.0  0.0001944965703569205  \\
            14.0  4.0354153436122295e-5  \\
            16.0  7.211403831424133e-6  \\
            18.0  9.434851449963251e-7  \\
            20.0  9.458496022713206e-7  \\
        }
        ;
    \addlegendentry {PROPT-$\htwo$}
\end{axis}
\end{tikzpicture} \\
  (a) $\hinf$ errors for \texttt{FOM-CONS} & (b) $\htwo$ errors for \texttt{FOM-CONS} \vspace{0.3cm} \\
 \begin{tikzpicture}[/tikz/background rectangle/.style={fill={rgb,1:red,1.0;green,1.0;blue,1.0}, draw opacity={1.0}}, show background rectangle]
\begin{axis}[point meta max={nan}, point meta min={nan}, legend cell align={left}, legend columns={3}, title={}, title style={at={{(0.5,1)}}, anchor={south}, font={{\fontsize{14 pt}{18.2 pt}\selectfont}}, color={rgb,1:red,0.0;green,0.0;blue,0.0}, draw opacity={1.0}, rotate={0.0}}, legend style={color={rgb,1:red,0.0;green,0.0;blue,0.0}, draw opacity={1.0}, line width={1}, solid, fill={rgb,1:red,1.0;green,1.0;blue,1.0}, fill opacity={1.0}, text opacity={1.0}, font={{\fontsize{8 pt}{10.4 pt}\selectfont}}, text={rgb,1:red,0.0;green,0.0;blue,0.0}, cells={anchor={center}}, at={(0.5, 1.02)}, anchor={south}}, axis background/.style={fill={rgb,1:red,1.0;green,1.0;blue,1.0}, opacity={1.0}}, anchor={north west}, xshift={1.0mm}, yshift={-1.0mm}, clip mode={individual}, width={0.45\textwidth}, height={0.25\textwidth}, scaled x ticks={false}, xlabel={\scriptsize{Reduced model order $r$}}, x tick style={color={rgb,1:red,0.0;green,0.0;blue,0.0}, opacity={1.0}}, x tick label style={color={rgb,1:red,0.0;green,0.0;blue,0.0}, opacity={1.0}, rotate={0}}, xlabel style={}, xmajorgrids={true}, xmin={1.46}, xmax={20.54}, xtick={{2.0,6.0,10.0,14.0,18.0}}, xticklabels={{$2$,$6$,$10$,$14$,$18$}}, xtick align={inside}, xticklabel style={font={{\fontsize{8 pt}{10.4 pt}\selectfont}}, color={rgb,1:red,0.0;green,0.0;blue,0.0}, draw opacity={1.0}, rotate={0.0}}, x grid style={color={rgb,1:red,0.0;green,0.0;blue,0.0}, draw opacity={0.1}, line width={0.5}, solid}, axis x line*={left}, x axis line style={color={rgb,1:red,0.0;green,0.0;blue,0.0}, draw opacity={1.0}, line width={1}, solid}, scaled y ticks={false}, ylabel={\scriptsize{$\mathcal{H}_\infty$ error}}, y tick style={color={rgb,1:red,0.0;green,0.0;blue,0.0}, opacity={1.0}}, y tick label style={color={rgb,1:red,0.0;green,0.0;blue,0.0}, opacity={1.0}, rotate={0}}, ylabel style={}, ymode={log}, log basis y={10}, ymajorgrids={true}, ymin={0.001410442397292485}, ymax={2.8757162348333445}, ytick={{1.0,0.1,0.01}}, yticklabels={{$10^{0}$,$10^{-1}$,$10^{-2}$}}, ytick align={inside}, yticklabel style={font={{\fontsize{8 pt}{10.4 pt}\selectfont}}, color={rgb,1:red,0.0;green,0.0;blue,0.0}, draw opacity={1.0}, rotate={0.0}}, y grid style={color={rgb,1:red,0.0;green,0.0;blue,0.0}, draw opacity={0.1}, line width={0.5}, solid}, axis y line*={left}, y axis line style={color={rgb,1:red,0.0;green,0.0;blue,0.0}, draw opacity={1.0}, line width={1}, solid}, colorbar={false}]
    \addplot[color={rgb,1:red,0.5961;green,0.3059;blue,0.6392}, name path={8858ae9f-99f9-4189-a7d1-e178f3cb7b63}, draw opacity={1.0}, line width={1}, solid, mark={star}, mark size={3.0 pt}, mark repeat={1}, mark options={color={rgb,1:red,0.0;green,0.0;blue,0.0}, draw opacity={1.0}, fill={rgb,1:red,0.5961;green,0.3059;blue,0.6392}, fill opacity={1.0}, line width={0.75}, rotate={0}, solid}]
        table[row sep={\\}]
        {
            \\
            2.0  1.191952716571951  \\
            4.0  0.6809347096819893  \\
            6.0  0.3888565852480281  \\
            8.0  0.09950281103759259  \\
            10.0  0.08878376435268408  \\
            12.0  0.08541062848364224  \\
            14.0  0.06924612876325249  \\
            16.0  0.016776027178113728  \\
            18.0  0.009885063869119502  \\
            20.0  0.005278262742707518  \\
        }
        ;
    \addplot[color={rgb,1:red,0.8941;green,0.102;blue,0.1098}, name path={387b2416-bb5e-4dae-b5f9-ceb785ef7b01}, draw opacity={1.0}, line width={1}, solid, mark={square*}, mark size={3.0 pt}, mark repeat={1}, mark options={color={rgb,1:red,0.0;green,0.0;blue,0.0}, draw opacity={1.0}, fill={rgb,1:red,0.8941;green,0.102;blue,0.1098}, fill opacity={1.0}, line width={0.75}, rotate={0}, solid}]
        table[row sep={\\}]
        {
            \\
            2.0  1.0875824001864745  \\
            4.0  0.6409  \\
            6.0  0.6262702929692878  \\
            8.0  0.16794945381752477  \\
            10.0  0.2226961391166911  \\
            12.0  0.15019286385882769  \\
            14.0  0.14486601238084534  \\
            16.0  0.1610328453692489  \\
            18.0  0.07747084563190718  \\
            20.0  0.015736188705940662  \\
        }
        ;
    \addplot[color={rgb,1:red,0.651;green,0.3373;blue,0.1569}, name path={79974278-0fdf-4913-9ae7-b112849ecd92}, draw opacity={1.0}, line width={1}, solid, mark={triangle*}, mark size={3.0 pt}, mark repeat={1}, mark options={color={rgb,1:red,0.0;green,0.0;blue,0.0}, draw opacity={1.0}, fill={rgb,1:red,0.651;green,0.3373;blue,0.1569}, fill opacity={1.0}, line width={0.75}, rotate={0}, solid}]
        table[row sep={\\}]
        {
            \\
            2.0  1.2025768518009847  \\
            4.0  0.6713557665068057  \\
            6.0  0.39144027714536056  \\
            8.0  0.1006522120458887  \\
            10.0  0.0894763636966087  \\
            12.0  0.08417227168678705  \\
            14.0  0.06727668627792127  \\
            16.0  0.01663938869600356  \\
            18.0  0.009824770456187057  \\
            20.0  0.005199949476551577  \\
        }
        ;
    \addplot[color={rgb,1:red,1.0;green,0.498;blue,0.0}, name path={2c5ec591-f49f-440d-92f5-a7a78b9997d0}, draw opacity={1.0}, line width={1}, solid, mark={triangle*}, mark size={3.0 pt}, mark repeat={1}, mark options={color={rgb,1:red,0.0;green,0.0;blue,0.0}, draw opacity={1.0}, fill={rgb,1:red,1.0;green,0.498;blue,0.0}, fill opacity={1.0}, line width={0.75}, rotate={180}, solid}]
        table[row sep={\\}]
        {
            \\
            2.0  2.31784307193479  \\
            4.0  1.8301514326834782  \\
            6.0  0.9565330491816848  \\
            8.0  0.6829899319057705  \\
            10.0  0.5309253654328733  \\
            12.0  0.3521299159121614  \\
            14.0  0.2699156915672576  \\
            16.0  0.19071171126410602  \\
            18.0  0.14741495043293967  \\
            20.0  0.08689288777193352  \\
        }
        ;
    \addplot[color={rgb,1:red,0.302;green,0.6863;blue,0.2902}, name path={90a1b374-331c-475d-9f34-f63c93f565ea}, draw opacity={1.0}, line width={1}, solid, mark={*}, mark size={3.0 pt}, mark repeat={1}, mark options={color={rgb,1:red,0.0;green,0.0;blue,0.0}, draw opacity={1.0}, fill={rgb,1:red,0.302;green,0.6863;blue,0.2902}, fill opacity={1.0}, line width={0.75}, rotate={0}, solid}]
        table[row sep={\\}]
        {
            \\
            2.0  0.6244657786948647  \\
            4.0  0.26723123777940055  \\
            6.0  0.19420046564853471  \\
            8.0  0.06277291424005145  \\
            10.0  0.04002650772668864  \\
            12.0  0.030555135431345696  \\
            14.0  0.021660427101826904  \\
            16.0  0.01035276257612093  \\
            18.0  0.0050554487344854375  \\
            20.0  0.0017499166139860977  \\
        }
        ;
    \addplot[color={rgb,1:red,0.2157;green,0.4941;blue,0.7216}, name path={30f60a4b-b481-4e5a-9cb9-b8144ab3a75a}, draw opacity={1.0}, line width={1}, solid, mark={diamond*}, mark size={3.0 pt}, mark repeat={1}, mark options={color={rgb,1:red,0.0;green,0.0;blue,0.0}, draw opacity={1.0}, fill={rgb,1:red,0.2157;green,0.4941;blue,0.7216}, fill opacity={1.0}, line width={0.75}, rotate={0}, solid}]
        table[row sep={\\}]
        {
            \\
            2.0  1.0768972987237624  \\
            4.0  0.6445  \\
            6.0  0.6514769563256232  \\
            8.0  0.17122649714091398  \\
            10.0  0.12807143278961117  \\
            12.0  0.15142350135415428  \\
            14.0  0.032119449405423044  \\
            16.0  0.01693458070754204  \\
            18.0  0.025132770938587807  \\
            20.0  0.01285244353641692  \\
        }
        ;
\end{axis}
\end{tikzpicture} & \begin{tikzpicture}[/tikz/background rectangle/.style={fill={rgb,1:red,1.0;green,1.0;blue,1.0}, draw opacity={1.0}}, show background rectangle]
\begin{axis}[point meta max={nan}, point meta min={nan}, legend cell align={left}, legend columns={3}, title={}, title style={at={{(0.5,1)}}, anchor={south}, font={{\fontsize{14 pt}{18.2 pt}\selectfont}}, color={rgb,1:red,0.0;green,0.0;blue,0.0}, draw opacity={1.0}, rotate={0.0}}, legend style={color={rgb,1:red,0.0;green,0.0;blue,0.0}, draw opacity={1.0}, line width={1}, solid, fill={rgb,1:red,1.0;green,1.0;blue,1.0}, fill opacity={1.0}, text opacity={1.0}, font={{\fontsize{8 pt}{10.4 pt}\selectfont}}, text={rgb,1:red,0.0;green,0.0;blue,0.0}, cells={anchor={center}}, at={(0.5, 1.02)}, anchor={south}}, axis background/.style={fill={rgb,1:red,1.0;green,1.0;blue,1.0}, opacity={1.0}}, anchor={north west}, xshift={1.0mm}, yshift={-1.0mm}, clip mode={individual}, width={0.45\textwidth}, height={0.25\textwidth}, scaled x ticks={false}, xlabel={\scriptsize{Reduced model order $r$}}, x tick style={color={rgb,1:red,0.0;green,0.0;blue,0.0}, opacity={1.0}}, x tick label style={color={rgb,1:red,0.0;green,0.0;blue,0.0}, opacity={1.0}, rotate={0}}, xlabel style={}, xmajorgrids={true}, xmin={1.46}, xmax={20.54}, xtick={{2.0,6.0,10.0,14.0,18.0}}, xticklabels={{$2$,$6$,$10$,$14$,$18$}}, xtick align={inside}, xticklabel style={font={{\fontsize{8 pt}{10.4 pt}\selectfont}}, color={rgb,1:red,0.0;green,0.0;blue,0.0}, draw opacity={1.0}, rotate={0.0}}, x grid style={color={rgb,1:red,0.0;green,0.0;blue,0.0}, draw opacity={0.1}, line width={0.5}, solid}, axis x line*={left}, x axis line style={color={rgb,1:red,0.0;green,0.0;blue,0.0}, draw opacity={1.0}, line width={1}, solid}, scaled y ticks={false}, ylabel={\scriptsize{$\mathcal{H}_2$ error}}, y tick style={color={rgb,1:red,0.0;green,0.0;blue,0.0}, opacity={1.0}}, y tick label style={color={rgb,1:red,0.0;green,0.0;blue,0.0}, opacity={1.0}, rotate={0}}, ylabel style={}, ymode={log}, log basis y={10}, ymajorgrids={true}, ymin={0.0013657660155785003}, ymax={3.834581055893011}, ytick={{1.0,0.1,0.01}}, yticklabels={{$10^{0}$,$10^{-1}$,$10^{-2}$}}, ytick align={inside}, yticklabel style={font={{\fontsize{8 pt}{10.4 pt}\selectfont}}, color={rgb,1:red,0.0;green,0.0;blue,0.0}, draw opacity={1.0}, rotate={0.0}}, y grid style={color={rgb,1:red,0.0;green,0.0;blue,0.0}, draw opacity={0.1}, line width={0.5}, solid}, axis y line*={left}, y axis line style={color={rgb,1:red,0.0;green,0.0;blue,0.0}, draw opacity={1.0}, line width={1}, solid}, colorbar={false}]
    \addplot[color={rgb,1:red,0.5961;green,0.3059;blue,0.6392}, name path={3a5c8562-727b-46b4-8833-8b57b5455208}, draw opacity={1.0}, line width={1}, solid, mark={star}, mark size={3.0 pt}, mark repeat={1}, mark options={color={rgb,1:red,0.0;green,0.0;blue,0.0}, draw opacity={1.0}, fill={rgb,1:red,0.5961;green,0.3059;blue,0.6392}, fill opacity={1.0}, line width={0.75}, rotate={0}, solid}]
        table[row sep={\\}]
        {
            \\
            2.0  1.1858688046711403  \\
            4.0  0.3437204865167503  \\
            6.0  0.2758572507802732  \\
            8.0  0.08485046536698404  \\
            10.0  0.0725959247092513  \\
            12.0  0.043057623849589687  \\
            14.0  0.031007866531573  \\
            16.0  0.016187996743688288  \\
            18.0  0.006461804613562718  \\
            20.0  0.0023035596091182564  \\
        }
        ;
    \addplot[color={rgb,1:red,0.8941;green,0.102;blue,0.1098}, name path={45bc6343-a470-4c9d-8615-1a69bb0bcfca}, draw opacity={1.0}, line width={1}, solid, mark={square*}, mark size={3.0 pt}, mark repeat={1}, mark options={color={rgb,1:red,0.0;green,0.0;blue,0.0}, draw opacity={1.0}, fill={rgb,1:red,0.8941;green,0.102;blue,0.1098}, fill opacity={1.0}, line width={0.75}, rotate={0}, solid}]
        table[row sep={\\}]
        {
            \\
            2.0  1.0820450865155542  \\
            4.0  0.29445975792017404  \\
            6.0  0.15054732290313463  \\
            8.0  0.08081838426350406  \\
            10.0  0.047900530219571184  \\
            12.0  0.02162302698223763  \\
            14.0  0.015285084003364462  \\
            16.0  0.009172313343568392  \\
            18.0  0.003207928943196011  \\
            20.0  0.00193261635283717  \\
        }
        ;
    \addplot[color={rgb,1:red,0.651;green,0.3373;blue,0.1569}, name path={4dff334a-eb18-4208-9b34-b526e5003ce6}, draw opacity={1.0}, line width={1}, solid, mark={triangle*}, mark size={3.0 pt}, mark repeat={1}, mark options={color={rgb,1:red,0.0;green,0.0;blue,0.0}, draw opacity={1.0}, fill={rgb,1:red,0.651;green,0.3373;blue,0.1569}, fill opacity={1.0}, line width={0.75}, rotate={0}, solid}]
        table[row sep={\\}]
        {
            \\
            2.0  1.2034777238970844  \\
            4.0  0.34600178534938014  \\
            6.0  0.27631197845823485  \\
            8.0  0.08541605722545421  \\
            10.0  0.07306280091512529  \\
            12.0  0.04438995514365575  \\
            14.0  0.030298092040109276  \\
            16.0  0.01605878337776877  \\
            18.0  0.006621881997463622  \\
            20.0  0.0023439485626462974  \\
        }
        ;
    \addplot[color={rgb,1:red,1.0;green,0.498;blue,0.0}, name path={f9c43d3e-b9af-4194-9d76-3b1f8209b0fe}, draw opacity={1.0}, line width={1}, solid, mark={triangle*}, mark size={3.0 pt}, mark repeat={1}, mark options={color={rgb,1:red,0.0;green,0.0;blue,0.0}, draw opacity={1.0}, fill={rgb,1:red,1.0;green,0.498;blue,0.0}, fill opacity={1.0}, line width={0.75}, rotate={180}, solid}]
        table[row sep={\\}]
        {
            \\
            2.0  3.0628332153254556  \\
            4.0  2.184599350417236  \\
            6.0  1.2345457889653655  \\
            8.0  0.8681480848998413  \\
            10.0  0.6623993328560581  \\
            12.0  0.42244190681823884  \\
            14.0  0.3354225711914794  \\
            16.0  0.21163525902328711  \\
            18.0  0.16040837956713272  \\
            20.0  0.10147370051386244  \\
        }
        ;
    \addplot[color={rgb,1:red,0.302;green,0.6863;blue,0.2902}, name path={ad7e03a6-e773-4c63-9ef2-6ce99a3f2da4}, draw opacity={1.0}, line width={1}, solid, mark={*}, mark size={3.0 pt}, mark repeat={1}, mark options={color={rgb,1:red,0.0;green,0.0;blue,0.0}, draw opacity={1.0}, fill={rgb,1:red,0.302;green,0.6863;blue,0.2902}, fill opacity={1.0}, line width={0.75}, rotate={0}, solid}]
        table[row sep={\\}]
        {
            \\
            2.0  1.0813745190444541  \\
            4.0  0.2943865312046926  \\
            6.0  0.14364786562663034  \\
            8.0  0.08397941165284865  \\
            10.0  0.04789919590128612  \\
            12.0  0.02162243752837442  \\
            14.0  0.015281822605858177  \\
            16.0  0.009931094204001697  \\
            18.0  0.0032065851782889917  \\
            20.0  0.00170990064490446  \\
        }
        ;
    \addplot[color={rgb,1:red,0.2157;green,0.4941;blue,0.7216}, name path={52ef187a-e2ac-4fbb-ae8c-c0fc9da0fedf}, draw opacity={1.0}, line width={1}, solid, mark={diamond*}, mark size={3.0 pt}, mark repeat={1}, mark options={color={rgb,1:red,0.0;green,0.0;blue,0.0}, draw opacity={1.0}, fill={rgb,1:red,0.2157;green,0.4941;blue,0.7216}, fill opacity={1.0}, line width={0.75}, rotate={0}, solid}]
        table[row sep={\\}]
        {
            \\
            2.0  1.0813686822537687  \\
            4.0  0.29438650285603046  \\
            6.0  0.14364783519834592  \\
            8.0  0.08078387074451146  \\
            10.0  0.058262991245602365  \\
            12.0  0.02162243168336662  \\
            14.0  0.02067025939699918  \\
            16.0  0.014098582692322463  \\
            18.0  0.005258975217113075  \\
            20.0  0.0019308357577836337  \\
        }
        ;
\end{axis}
\end{tikzpicture} \\
  (c) $\hinf$ errors for \texttt{FOM-RAND} & (d) $\htwo$ errors for \texttt{FOM-RAND}
    \end{tabular}
    \caption{$\hinf$ and $\htwo$ error comparison for \texttt{FOM-CONS} and \texttt{FOM-RAND}\@. Note that in the $\hinf$ error comparisons SOBMOR-$\hinf$ and in the $\htwo$ error comparisons SOBMOR-$\htwo$ is used.}\label{fig:FOMRANDORIGComp}
\end{figure*}

The $\htwo$ errors, reported in Figures~\ref{fig:FOMRANDORIGComp} (b) and (d), exhibit a less distinct behavior. Again, it can be clearly seen that IRKA-PH has by far the worst accuracy. The other methods have similar and much better accuracies --- in particular for the \emph{simpler} system \texttt{FOM-CONS}\@. For \texttt{FOM-RAND}, the $\hinf$ methods XminBT and PRBT (which once again have similar errors), have a slightly worse performance for most reduced model orders. Note that in this comparison, we use SOBMOR-$\htwo$ instead of SOBMOR-$\hinf$ because the $\htwo$ errors of the models obtained using SOBMOR-$\hinf$ are infinite because the feedthrough terms are not matched.

\begin{figure}[htpb]
  \centering
    \begin{tabular}{c}
      \input{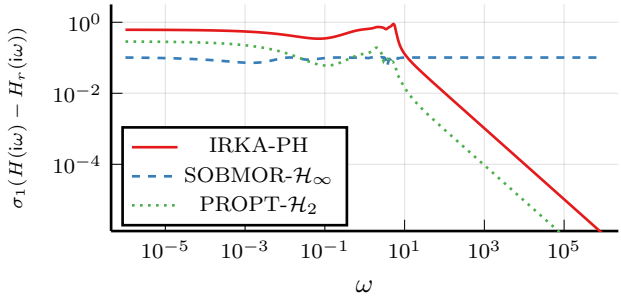} \\
      \phantom{(aaaaa)} (a) Errors for ROM dimension $r=10$ \\
      \phantom{(a)} \\
      \input{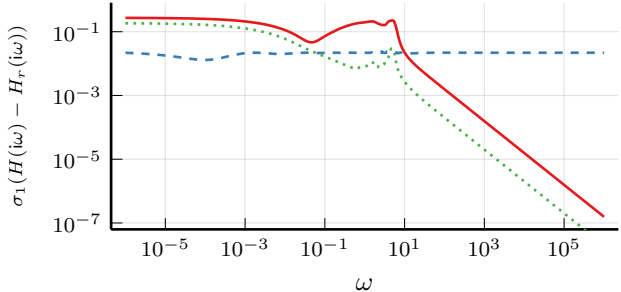} \\
      \phantom{(aaaaa)} (b) Errors for ROM dimension $r=20$
    \end{tabular}
  \caption{Error sigma plots for \texttt{FOM-MIMO}}\label{fig:FOM_XXL_Freq}
\end{figure}

In Figure~\ref{fig:FOM_XXL_Freq}, we report the error transfer functions between \texttt{FOM-MIMO} and the ROMs obtained with IRKA-PH and our proposed methods. Due to the vast system dimension of \texttt{FOM-MIMO}, we do not apply the other methods because no implementations of these methods that exploit sparsity are currently publicly available. Furthermore, the exact computation of $\hinf$ or $\htwo$ errors is computationally prohibitive. However, the error transfer functions indicate that our methods continue to work as intended even in the large-scale case. In particular SOBMOR-$\hinf$ leads to a flat error curve in the sigma plot, which has its highest peak value well-below the other errors and PROPT-$\htwo$ has an error transfer function that is below the error of IRKA-PH over the entire imaginary axis and below the error of SOBMOR-$\hinf$ for higher frequencies.

\begin{rem}
  Note that the passivity-preserving methods that we use for comparison with our methods require a transformation as in Theorem~\ref{thm:KYP}, in order to recover the pH structure from the computed ROM\@. In all of our experiments, numerical inaccuracies lead to a passivity matrix of the transformed system that is not positive semi-definite but rather has a few slightly negative eigenvalues. These are typically of the order of machine precision. However, in some cases, the smallest negative eigenvalues have absolute values of the order of $10^{-11}$. Only IRKA-PH and our optimization-based methods ensure a pH structure with a positive semi-definite passivity matrix.
\end{rem}

\section{Conclusion}%
\label{sec:conclusion}

We have presented two optimization-based methods for structure-preserving MOR of pH-DAEs. These make it possible to compute accurate ROMs with respect to either the $\hinf$ or the $\htwo$ norm. The main benefits compared to state-of-the-art methods are the simplified treatment of the algebraic equations, which can be incorporated into the parameterized ROM in a structure-preserving way without increasing its state dimension. Furthermore, our methods are data-driven, such that no transformations to the FOM system matrices are required. Nonetheless, we have shown, how transformations can be applied in order to obtain an accurate estimate of the feedthrough, which is essential in the $\htwo$ case. Finally, our numerical experiments show that the optimization-based methods often lead to a higher accuracy (especially in the $\hinf$ norm).

We are currently investigating the application of our method to higher index pH-DAEs. These may have improper transfer functions, which are not currently supported in our parameterization.

\bibliographystyle{elsarticle-num} 
\bibliography{references}

\section*{Appendix}
In Algorithm~\ref{alg:sobmor_trapz} we give details for the integral computation of~\eqref{eq:htwonorm}. It implements an adaptive trapezoidal rule, in which new quadrature points are added at the logarithmic midpoint of a given interval if the relative accuracy requirement of the integral over the given interval is not met. In our implementation, we initialize the interval stack with a list of intervals that has been used in the previous function call and cache the function evaluations in line 6 for subsequent iterations. In our implementation of SOBMOR-$\htwo$, we integrate the function
\begin{align}
  \label{eq:htwointegral}
  f : \R  \rightarrow \R, \qquad \omega \mapsto {\| \pHtf(\ri \omega) - \pHtfr(\ri \omega, \theta) \|}_{\rm F}^2
\end{align}
over the interval $[0, 10^8]$ and use the integral as an objective function for the minimization of the $\htwo$ error.

\begin{algorithm}[tbh]
  \LinesNumbered
  \SetAlgoLined
  \DontPrintSemicolon
  \SetKwInOut{Input}{Input}\SetKwInOut{Output}{Output}
  \caption{Adaptive Quadrature}\label{alg:sobmor_trapz}
  \Input{Function $f : \R  \rightarrow \R$ in \eqref{eq:htwointegral},
    initial interval stack $S_I$, relative error tolerance $\epsilon_I > 0$.
  }
  \Output{Approximate integral $I$ of $f$ over all intervals in $S_I$.}
  Initialize $I := 0$. \;
  \While{$S_I$ is not empty}{
    Pop first interval $I$ from $S_I$. \;
    Set $\alpha, \beta$ with $\alpha < \beta$ as endpoints of $I$. \;
    Set $\gamma := \exp(\ln(10) \log_{10}((\alpha+\beta)/2))$. \;
    Set $f_\alpha := f(\alpha),\, f_\beta = f(\beta), \,f_\gamma = f(\gamma)$. \;
    Set $S_1 := |\alpha-\beta|(f_\alpha + f_\beta)/2$.\;
    Set $S_2 := (|\alpha-\gamma|(f_\alpha + f_\gamma)+|\beta-\gamma|(f_\beta+f_\gamma))/2$.\;
    \eIf{$|S_2-S_1|/|S_2| < \epsilon_I$}{
      Set $I := I + S_2$.
      }{
      Push intervals $[\alpha, \gamma]$ and $[\gamma, \beta]$ to $S_I$. \;
    }
  }
\end{algorithm}

\end{document}